\documentclass{article}

\usepackage{arxiv}

\usepackage[utf8]{inputenc} 
\usepackage[T1]{fontenc}    
\usepackage{hyperref}       
\usepackage{url}            
\usepackage{booktabs}       
\usepackage{amsfonts}       
\usepackage{nicefrac}       
\usepackage{microtype}      
\usepackage{doi}

\usepackage{graphicx}
\usepackage{amsmath}
\usepackage{amssymb}
\usepackage{enumitem}
\usepackage{float}
\usepackage{booktabs}
\usepackage{xcolor}
\usepackage{subcaption}
\usepackage{bm}
\usepackage{textcomp}
\usepackage{geometry}
\usepackage{multirow}
\usepackage{makecell}
\usepackage{booktabs}

\usepackage{siunitx}
\usepackage[T1]{fontenc}

\DeclareMathOperator{\sgn}{sgn}

\usepackage{makecell}

\usepackage{amsthm}

\theoremstyle{plain}
\newtheorem{thm}{Theorem}[section]

\newtheorem{prop}[thm]{Proposition}
\newtheorem*{cor}{Corollary}
\theoremstyle{definition}
\newtheorem{dfn}{Definition}

\theoremstyle{remark}
\newtheorem*{remark}{Remark}

\newtheorem{asm}{Assumption}

\title{Calculation of explicit expressions for the Hopf bifurcation limit cycles in delay-differential equations}

\author{ \href{https://orcid.org/0009-0002-0496-1282}{\includegraphics[scale=0.06]{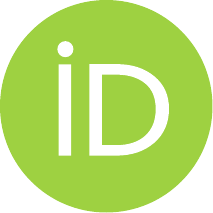}\hspace{1mm}José ~Enríquez Gabeiras}\\
	Department of Applied Mathematics, Universidad
	Politécnica de Madrid.\\
	Madrid, Spain \\
	\texttt{jose.enriquez@alumnos.upm.es} \\
	\And
	\href{https://orcid.org/0000-0001-7092-3554}{\includegraphics[scale=0.06]{orcid.pdf}\hspace{1mm}Juan Francisco Padial Molina} \\
	Department of Applied Mathematics, Universidad
	Politécnica de Madrid.\\
	Madrid, Spain \\
	\texttt{jfpadial@upm.es} \\
}

\date{}

\hypersetup{
	pdftitle={Calculation of explicit expressions for the Hopf bifurcation limit cycles in delay-differential equations},
	pdfsubject={q-bio.NC, q-bio.QM},
	pdfauthor={José ~Enríquez Gabeiras, Juan F. Padial Molina},
	pdfkeywords={delay-differential equations, Hopf bifurcation, Poincaré-Lindstedt method, car-following model, SIR epidemic model},
}

\begin{document}
	
	\maketitle
	
	\begin{abstract}
		This paper introduces a methodology to derive explicit power series approximations for the limit cycle periodic solutions of the Hopf bifurcation in autonomous discrete delay differential equations (\textit{DDE}). The procedure extends the methodology introduced in \cite{Casal1980} (\textit{Casal $ \& $ Freedman}, 1980) by providing a detailed algorithm that iteratively performs systematic calculations up to any desired order of approximation, ensuring a specific error tolerance for any nonlinear \textit{DDE}  presenting a Hopf bifurcation. The methodology is applied to two relevant delay-differential models to illustrate its features: a recently introduced car-following mobility model, whose oscillations are a plausible explanation for the density waves and congestion in road traffic, and a SIR epidemic model for propagation of diseases with temporary immunity.
	\end{abstract}

	\keywords{delay-differential equations, Hopf bifurcation, Poincaré-Lindstedt method, car-following model, SIR epidemic model}

\section{Introduction}
We consider the autonomous nonlinear delay differential equation for the vector function $ \mathbf{x} \in \mathbb{R}^n $,  $ n \geq 1 $:
\begin{align}\label{ecu1}
	\mathbf{x}'(t) = \mathbf{g}\left( \lambda, \mathbf{x}(t), \mathbf{x}(t-\lambda)\right), 
\end{align}
where $ t \in \mathbb{R}^{+} $, and $ \lambda \in \mathbb{R}^{+} $ is some parameter. We also consider $ \mathbf{g} \in C^{k}(\mathbb{R}^{+} \times \mathbb{R}^n \times \mathbb{R}^n;\mathbb{R}^n) $ with $ k \geq 3 $, and, without loss of generality, the equilibrium point of the system to be $ \mathbf{x}(t) = \mathbf{0} $, at which $ \mathbf{g}(\lambda, \mathbf{0}, \mathbf{0}) = \mathbf{0} \hspace{0.3 cm} \forall \lambda \geq 0 $.

This equation is a particular case of the \textit{Retarded Functional Differential Equation (RFDE)} class defined in \cite{Hale1977}, where, in this instance, $ \textbf{g} $ presents a single, constant delay $ \lambda $. These equations are also known as constant discrete delay equations \cite{Smith2011}, and will be subsequently referred to as \textit{DDE}. In \cite{Hale1977} conditions for existence, uniqueness, continuity, differentiability and continuation of the equation solutions are provided. The solutions to this type of equations can undergo the Hopf bifurcation phenomenon when the parameter $ \lambda $ is varied, such that starting from a specific threshold $ \lambda_{0} $ the solution changes from a constant equilibrium to a nonzero periodic solution. In this paper we provide a method to systematically derive explicit expressions for the oscillatory solutions of the Hopf bifurcation of equations of type \eqref{ecu1} up to any desired approximation order, based on the original method introduced in \cite{Casal1980}.

In order to assess the stability of the system, the linearized approximation of \eqref{ecu1} is used:
\begin{align}\label{ecu2} 
	\pmb{\eta}'(t) = \mathbf{P}(\lambda) \pmb{\eta}(t) + \mathbf{Q}(\lambda) \pmb{\eta}(t-\lambda), 
\end{align}
where $ \mathbf{P}(\lambda) = \mathbb{J}_{\mathbf{u}} \mathbf{g}(\lambda, \mathbf{u}, \mathbf{v}) \mid_{\mathbf{u}=\mathbf{v}=\mathbf{0}} $ and $ \mathbf{Q}(\lambda) = \mathbb{J}_{\mathbf{v}} \mathbf{g}(\lambda, \mathbf{u}, \mathbf{v}) \mid_{\mathbf{u}=\mathbf{v}=\mathbf{0}} $, and $ \mathbb{J}_{*} $ is the jacobian matrix with respect to the vector $ * $. Based on \cite{Hale1977, Smith2011}, the following assumptions are considered:
\begin{asm}[A1]\label{A1} The linear system \eqref{ecu2} has a non-constant periodic solution for some $ \lambda = \lambda_{0} $, which is determined by a pair of purely imaginary roots $ \left\lbrace r_{0}= j\omega_{0} , \overline{r}_{0}= -j\omega_{0}\right\rbrace , \omega_{0} > 0 $, in the characteristic equation of \eqref{ecu2}, while the other roots of the characteristic equation $ r_{i} \neq r_{0}, \overline{r_{0}} $ satisfy $ r_{i} \neq m r_{0} $ for any integer $ m $.
\end{asm}
\begin{asm}[A2]\label{A2} By the differentiability of $ \mathbf{g} $ and the implicit function theorem, there exist a continuously differentiable family of roots $ r(\lambda) $ for small deviations of $ \lambda $ for which $ Re(r'(\lambda)) \neq 0 $, implying that the roots $ r(\lambda) $ cross the imaginary axis as $ \lambda $ varies.
\end{asm}
These assumptions imply that for small deviations of $ \lambda_{0} $ there exist non-constant periodic solutions for the nonlinear equation \eqref{ecu1}  with period close to $ 2 \pi / \omega_{0} $, as stated by the \textit{Hopf Bifurcation Theorem} (see \textit{Theorem 1.1} in ch. 11 of \cite{Hale1977}, or \textit{Theorem 6.1} in \cite{Smith2011}).

Methods to derive numerical solutions for this type of bifurcations have been widely used for a long time. In particular, owing to the mathematical tractability of the singular perturbation methods like multiple scales, harmonic balance or Poincaré-Lindstedt, they have traditionally been the choice to find approximate explicit solutions to ODEs bifurcations \cite{Verhulst1990,Nayfeh2008}. More recently they have also been adapted to be applied to other type of equations, like \textit{DDEs} and fractional calculus equations \cite{ Min2010_CiteCF80, Rand2011, Zhongjin2014_CiteCF80, Min_2015_CiteCF80, Bridgewater2020_CiteCF80}. Due to the complexity of the calculation of the higher order terms in these methods, the majority of applications have been restricted to systems of equations derived from an scalar differential equation (comparatively easier to manipulate), and the validity of the solution  has been only limited to small deviations from the limit value of the bifurcation parameter.

The reference \cite{Casal1980} established a new framework to generically tackle \textit{DDE} bifurcations by extending the Poincaré-Lindstedt method for ODEs to  any nonlinear autonomous \textit{DDE} equation that complies with specific solvability conditions. The paper provides an existence theorem for uniformly valid asymptotic series expansions of the solution, alongside a detailed recursive algorithm for the calculation of such expansions. However, no methodical approach to solve the \textit{DDE} that arise in each iteration step was provided, thus limiting the applicability of the method to lower order systems up to few orders of approximation, where these calculations are manageable. This paper extends that methodology with an algorithm to systematically perform the calculations of these intermediate \textit{DDEs}, thus yielding approximate periodic solutions for equations like \eqref{ecu1} up to any desired expansion order, allowing the calculation of accurate approximations even for delays significantly greater than the bifurcation delay $ \lambda_{0} $.

To illustrate the complete procedure enabled by this extended methodology, we present explicit expressions for the approximate periodic solutions for two relevant problems in the area of delayed systems. First, the original analysis of a new car-following model recently introduced in \cite{Padial2022} (\textit{nDDE} model) is expanded by providing explicit expressions for its periodic equilibrium solutions. Secondly, we include results for an epidemic disease transmission model (SIR model), in which temporary immunity induces oscillations and subsequent epidemic waves, and for which an explicit relationship between the immunity period and the frequency and strength of these waves can be derived using this methodology.

\section{Poincaré-Lindstedt method applied to \textit{DDEs}}\label{PL_method}

\subsection{System of equations transformation}\label{sys_eq}
 The Poincaré-Lindstedt method for \textit{DDEs} requires a number of sequential transformations to the original equation so that any \textit{DDE} problem is expressed in a general standard form. First we transform equations \eqref{ecu1} and \eqref{ecu2} by means of a re-scaling of the time variable, $ \hat t = t \omega_{0} $, where $ \omega_{0} $ is the fundamental frequency of the Hopf bifurcation (see Assumption \ref{A1}). Thus we get the $ DDE_{\hat t} $ form of the problem:
\begin{align}\label{ecu3}
	\mathbf{x}'(\hat t) = \frac{\mathbf{g}(\hat \lambda, \mathbf{x}(\hat t), \mathbf{x}(\hat t - \hat \lambda))}{\omega_{0}} = \mathbf{f}(\hat \lambda, \mathbf{x}(\hat t), \mathbf{x}(\hat t - \hat \lambda)), 
\end{align}
\begin{align}\label{ecu4} 
	\pmb{\eta}'(\hat t) = \mathbf{A}(\hat \lambda) \pmb{\eta}(\hat t) + \mathbf{B}(\hat \lambda) \pmb{\eta}(\hat t-\hat \lambda), 
\end{align}
where $ \hat \lambda = \lambda\omega_{0} $, with $ \mathbf{A}(\hat \lambda)=\mathbf{P}(\lambda) / \omega_{0} $ and $ \mathbf{B}(\hat \lambda)=\mathbf{Q}(\lambda) / \omega_{0} $.
The particularization of \eqref{ecu4} for $ \hat{\lambda} = \hat{\lambda_{0}} = \lambda_{0} \omega_{0} $ has periodic solutions with period $ 2 \pi $, among which we will consider a solution $ \mathbf{z}_{0}(\hat t) $ that verifies:
\begin{align}\label{ecu5} 
	\dfrac{1}{2 \pi}\int_{0}^{2 \pi}|\mathbf{z}_{0}(\hat t) |^{2} d \hat t = 1.
\end{align}
Additionally, we denote by $ L $ the constant coefficient differential delay operator derived from \eqref{ecu4} for $ \hat{\lambda} = \hat{\lambda_{0}} $, that will be used in the iterative calculations in later steps of the process:
\begin{align}\label{ecu6} 
	L \pmb{\eta} = \pmb{\eta}'(\hat t) - \mathbf{A}(\hat \lambda_{0}) \pmb{\eta}(\hat t) - \mathbf{B}(\hat \lambda_{0}) \pmb{\eta}(\hat t - \hat \lambda_{0}), 
\end{align}
 and we denote the null space of the operator $ L $ as $ N(L) \subset \mathbb{P}^{\infty}$, ($ \mathbb{P}^{\infty} $ being the subspace of $ C^{\infty} $ of $ 2 \pi $ periodic vector functions). In a general bifurcation situation, $ N(L) $ could have a given dimension $ l \geq 2 $, determined by the nature of the bifurcation. In the specific case of the Hopf bifurcation, the single pair of eigenvalues assumed in Assumption \ref{A1} determines a dimension $ l = 2 $, so we will consider a fixed orthonormal basis $ \{\mathbf{v}_{1}(\hat t), \mathbf{v}_{2}(\hat t)\} $ for the space $ N(L) $.
 
\noindent It will also be necessary to consider the formal adjoint operator $ L^{*} $:
\begin{align}\label{ecu7}
	L^{*} \pmb{\eta} = \pmb{\eta}'(\hat t) + \mathbf{A}^{T}(\hat \lambda_{0}) \pmb{\eta}(\hat t) + \mathbf{B}^{T}(\hat \lambda_{0}) \pmb{\eta}(\hat t + \hat \lambda_{0}), 
\end{align}
with $ \dim(N(L^{*})) = \dim(N(L)) = 2 $ and a fixed orthonormal basis of the space $ N(L^{*}) $ $ \{\mathbf{w}_{1}(\hat t), \mathbf{w}_{2}(\hat t)\} $, $ \mathbf{A}^{T} $ and $ \mathbf{B}^{T} $ being the transposes of $ \mathbf{A} $ and $ \mathbf{B} $ respectively.

The Poincaré-Lindstedt method (either for ODEs or \textit{DDEs}) is based on the search of the periodic solutions of \eqref{ecu3} by means of an expansion parameter $ \varepsilon $, defined, for convenience, as:
\begin{align}\label{ecu8} 
	\dfrac{1}{\hat{T}}\int_{\hat{T}}|\mathbf{x}(\hat t) |^{2} d \hat t = \kappa (\varepsilon) \cdot \varepsilon^{2}, 
\end{align}
where $ \mathbf{x}(\hat t) $ is the periodic solution with period $ \hat{T} $ that arise in the Hopf bifurcation for the delay $ \hat \lambda $, and $ \kappa (\varepsilon) $ is a scaling factor that will be determined during the solution construction (see Remark \ref{rmk1}). In this way, for any $ \varepsilon > 0 $ we can change the dependent variable $ \mathbf{x}(\hat t) = \varepsilon \mathbf{z}(\hat t) $, so that
\eqref{ecu3} and \eqref{ecu8} are transformed into the $ DDE_{z} $ form of the problem :
\begin{align}\label{ecu9}
	\mathbf{z}'(\hat t)=
	\begin{cases}
		\dfrac{\mathbf{f}\left( \hat \lambda, \varepsilon \mathbf{z}(\hat t), \varepsilon \mathbf{z}(\hat t - \hat \lambda)\right) }{\varepsilon} & \text{ for 	} \varepsilon \neq 0, \\
		\\
		\mathbf{A}(\hat \lambda) \mathbf{z}(\hat t) + \mathbf{B}(\hat \lambda) \mathbf{z}(\hat t-\hat \lambda) & \text{ for } \varepsilon = 0, 
	\end{cases}
\end{align}
\begin{align}
	\nonumber \dfrac{1}{\hat{T}(\varepsilon)}\int_{0}^{\hat{T}(\varepsilon)}|\mathbf{z}(\hat t) |^{2} d \hat t = \kappa (\varepsilon), 
\end{align}
where the implicit dependence of the period $ \hat{T} $ on $ \varepsilon $ is noted.

Finally, the time variable is again re-scaled as $ \tau = (2 \pi / \hat{T}(\varepsilon)) \cdot \hat t $, so that the solution is periodic $ 2 \pi $ in $ \tau $ for any $ \hat \lambda(\varepsilon) > \hat \lambda_{0} $, noting now the implicit dependence of $ \hat \lambda $ on $ \varepsilon $. We then define:
\begin{align}
\nonumber \mathbf{Z}(\tau, \varepsilon) = \mathbf{z}\left( \frac{\hat{T}(\varepsilon)}{2 \pi}\tau\right), 
\end{align}
so, $ \mathbf{Z}(\tau, \varepsilon) \in \mathbb{P} $ (space of $ 2 \pi $ periodic vector functions). Moreover, we get the expression for the $ DDE_{Z} $ form of the problem:
\begin{align}\label{ecu12}
	\mathbf{Z}'(\tau, \varepsilon)=
	\begin{cases}
		\dfrac{\hat{T}(\varepsilon)}{2 \pi}\dfrac{\mathbf{f}\left( \hat \lambda(\varepsilon), \varepsilon \mathbf{Z}(\tau, \varepsilon), \varepsilon \mathbf{Z}\left( \tau - \frac{2 \pi}{\hat{T}(\varepsilon)} \hat \lambda(\varepsilon), \varepsilon\right) \right) }{\varepsilon} & \text{ for } \varepsilon \neq 0, \\
		\\
		\mathbf{A}(\hat \lambda_{0}) \mathbf{Z}(\tau, 0) + \mathbf{B}(\hat \lambda_{0}) \mathbf{Z}(\tau-\hat \lambda_{0}, 0) & \text{ for } \varepsilon = 0, 
	\end{cases}
\end{align}
\begin{align}\label{ecu13}
	\dfrac{1}{2 \pi}\int_{0}^{2 \pi}|\mathbf{Z}(\tau, \varepsilon) |^{2} d \tau = \kappa (\varepsilon).
\end{align}
Notice that for $ \varepsilon = 0 $ we have $ \mathbf{Z}(\tau, 0) = \mathbf{z}_0(\tau) $, and 
\begin{align}
	\nonumber \frac{1}{2 \pi} \int_{0}^{2 \pi}|\mathbf{Z}(\tau, 0)|^{2} d \tau = \frac{1}{2 \pi} \int_{0}^{2 \pi}|\mathbf{z}_{0}(\tau)|^{2} d \tau = 1 = \kappa(0).
\end{align}

 We finally select any linear functional $ \Phi: \mathbb{P} \longrightarrow \mathbb{R} $ that verifies:
\begin{align}
	\nonumber \Phi \mathbf{z}_{0} = 0, 
\end{align}
which allows us to define the additional equation:
\begin{align}\label{ecu16}
	\Phi \mathbf{Z} = 0.
\end{align}
The role of this functional is to set a specific time reference for the solution, given that we will build a periodic solution lacking any particular initial condition. This completes the system of equations \eqref{ecu12}, \eqref{ecu13} and \eqref{ecu16}, which is verified for $ \varepsilon = 0 $ by $ \mathbf{z}_{0}(\tau) $. 

\subsection{Solutions as series expansion}\label{sol_ser_exp}
Having set the original equations in terms of the new variables, we express the solutions of $ \mathbf{Z}(\tau, \varepsilon), \hat \lambda(\varepsilon) $ and $ \hat{T}(\varepsilon) $ as asymptotic power series expansions on $ \varepsilon $:
\begin{align}\label{ecu17} 
	\mathbf{Z}(\tau, \varepsilon) = \sum_{j=0}^{m} \mathbf{Z}_{j}(\tau)\varepsilon^{j} + O(\varepsilon^{m+1}), 
\end{align}
\begin{align}\label{ecu18} 
	 \hat \lambda(\varepsilon) = \sum_{j=0}^{m} \hat \lambda_{j}\varepsilon^{j} + O(\varepsilon^{m+1}), 
\end{align}
\begin{align}\label{ecu19} 
	\hat{T}(\varepsilon) = \sum_{j=0}^{m} \hat{T}_{j}\varepsilon^{j} + O(\varepsilon^{m+1}), 
\end{align}
where $ \mathbf{Z}_{j}: \mathbb{R}^{+} \rightarrow \mathbb{R}^n $ and $ \mathbf{Z}_{j}(\tau) \in \mathbb{P} $.
The values $ \hat \lambda_{0} = \lambda_{0} \cdot \omega_{0} $ and $ \hat{T}_{0} = 2 \pi $ are known. In addition we define the linear functional $ \Phi $ in \eqref{ecu16} as the one that zeroes the first component of the solution vector $ \mathbf{Z} $ in $ \tau = 0 $:
\begin{align}\label{ecu20}
	& \Phi \mathbf{Z} = \Phi \mathbf{Z}(\cdot, \varepsilon) = \mathbf{Z}^{1}(0, \varepsilon) = 0, 
\end{align}
 where the superscript $ []^{1} $ denotes the first component of a vector.
The calculation starts by determining $ \mathbf{Z}_{0}(\tau) = \mathbf{z}_{0}(\tau) $ from the basis $ \{\mathbf{v}_{1}, \mathbf{v}_{2}\} $ such that it verifies \eqref{ecu20} :
\begin{align}\label{ecu22}
	\mathbf{Z}_{0}^{1}(0) = 0.
\end{align}

In the following we summarize the methodology of iterative calculations introduced in \cite{Casal1980}.

\subsubsection{Equation for order $ j $ of $ \varepsilon $}

In order to set the notation, we first note the general expression of the Taylor expansion of $ \mathbf{f}(\hat \lambda, \mathbf{x}, \mathbf{y}) $ around the equilibrium point $ (\mathbf{x}_{0}, \mathbf{y}_{0}) = (\mathbf{0}, \mathbf{0}) $ as:
\begin{align}\label{ecu23} 
	\mathbf{f}(\hat \lambda, \mathbf{x}, \mathbf{y}) = \mathbf{A}(\hat \lambda) \mathbf{x} + \mathbf{B}(\hat \lambda) \mathbf{y} + C(\hat \lambda, \mathbf{x}, \mathbf{x}) + D(\hat \lambda, \mathbf{x}, \mathbf{y}) + E(\hat \lambda, \mathbf{y}, \mathbf{y}) + r(\hat \lambda, \mathbf{x}, \mathbf{y}), 
\end{align}
with $ C, D $ and $ E $ representing the quadratic terms of the expansion (bilinear forms), while $ r $ represents the remaining higher order terms.

The solution of the system is found by first introducing the series expansion of the variables (\eqref{ecu17} - \eqref{ecu19}) in the $ DDE_{Z} $ equation \eqref{ecu12}, and equating the terms for the same powers of $\varepsilon$. Then we get for $ j \geq 1 $:
\begin{align}\label{ecu24}
	\mathbf{Z}'_{j}  (\tau)= 
	\text{$ \varepsilon^j $ term of \textit{Series}} 
	\left( \frac{\hat{T}(\varepsilon)}{2 \pi}\frac{\mathbf{f}\left( \hat \lambda(\varepsilon), \varepsilon \mathbf{Z}(\tau, \varepsilon), \varepsilon \mathbf{Z}\left( \tau - \frac{2 \pi}{\hat{T}(\varepsilon)} \hat \lambda(\varepsilon), \varepsilon\right) \right) }{\varepsilon} \right), 
\end{align}
which, inserting in \eqref{ecu24} the expression of $ \mathbf{f} $ in \eqref{ecu23}, can be expressed for any $ j $ as the following nonhomogeneous linear \textit{DDE}, that we will call $ DDE_{Z_{j}} $:
\begin{align}\label{ecu25}
\mathbf{Z}_{j}' (\tau)=
		\mathbf{A}(\hat \lambda_{0}) \mathbf{Z}_{j}(\tau) + \mathbf{B}(\hat \lambda_{0}) \mathbf{Z}_{j}(\tau-\hat \lambda_{0}) + \mathbf{h}_{j}(\tau), 
\end{align}
where $ \mathbf{h}_{j}(\tau) $ depends on the terms $ \mathbf{Z}_{k}(\tau), k<j $, and verifies, by \eqref{ecu6}: 	
\begin{align}
	\nonumber \mathbf{h}_{j}(\tau) = L\mathbf{Z}_{j} \hspace{0.1 cm}.
\end{align}	
In order to express the structure of $ \mathbf{h}_{j}(\tau) $, we also define the operators $ M, N : \mathbf{V} \in \mathbb{R}^n \rightarrow \mathbb{R}^n $:
\begin{align}\label{ecu27} 
	M\mathbf{V} = C(\hat \lambda_{0}, \mathbf{Z}_{0}(\tau), \mathbf{V}) +\frac{ D(\hat \lambda_{0}, \mathbf{V}, \mathbf{Z}_{0}(\tau - \hat \lambda_{0}))}{2}, 
\end{align}
\begin{align}\label{ecu28} 
	N\mathbf{V} = \frac{ D(\hat \lambda_{0}, \mathbf{Z}_{0}(\tau ), \mathbf{V})}{2} + E(\hat \lambda_{0}, \mathbf{Z}_{0}(\tau - \hat \lambda_{0}), \mathbf{V}).
\end{align}	
In terms of the defined notation, the following relationship arises:
\begin{align}\label{ecu29} 
	\mathbf{h}_{j}(\tau) = L\mathbf{Z}_{j}(\tau) = \hat{T}_{j} \mathbf{R}(\tau) + \hat \lambda_{j} \mathbf{S}(\tau) + M \mathbf{Z}_{j-1}(\tau) + N \mathbf{Z}_{j-1}(\tau-\hat \lambda_{0}) + \mathbf{P}_{j}(\tau), 
\end{align}
with:
\begin{align}\label{ecu30} 
	\mathbf{R}(\tau) = \frac{1}{2 \pi}\left( \mathbf{Z}'_{0}(\tau) + \hat \lambda_{0} \mathbf{B}(\hat \lambda_{0}) \mathbf{Z}'_{0}(\tau -\hat \lambda_{0})\right), 
\end{align}
\begin{align}\label{ecu31} 
	\mathbf{S}(\tau) = \mathbf{A}'(\hat \lambda_{0}) \mathbf{Z}_{0}(\tau) + \mathbf{B}'(\hat \lambda_{0}) \mathbf{Z}_{0}(\tau-\hat \lambda_{0}) - \mathbf{B}(\hat \lambda_{0}) \mathbf{Z}'_{0}(\tau-\hat \lambda_{0}), 
\end{align}
where $ \mathbf{A}'(\hat \lambda)=\frac{\partial \mathbf{A}(\hat \lambda)}{\partial \hat \lambda} $, $ \mathbf{B}'(\hat \lambda)=\frac{\partial \mathbf{B}(\hat \lambda)}{\partial \hat \lambda} $, and $ \mathbf{P}_{j}(\tau) $ is a $ 2 \pi $ periodic function depending on the following variables:
\begin{align}
	\nonumber \mathbf{P}_{j}(\tau) \equiv \mathbf{P}_{j}(\mathbf{Z}_{0}(\tau), ..., \mathbf{Z}_{j-1}(\tau), \hat \lambda_{0}, ..., \hat \lambda_{j-1}, \hat{T}_{0}, ..., \hat{T}_{j-1}), 
\end{align}
with $ \mathbf{P}_{1} = \mathbf{0} $. So we get from \eqref{ecu29} a different \textit{DDE} for each order of $ \varepsilon $, which expresses $ \mathbf{Z}_{j} $ in terms of $ \hat \lambda_{k} $, $ \hat{T}_{k} $, and $ \mathbf{Z}_{k}(\tau), k=0..(j-1) $.

To calculate $ \mathbf{Z}_j(\tau) $ we use two additional equations. From \eqref{ecu20}:	  
\begin{align}\label{ecu33}
	& \mathbf{Z}_{j}^{1}(0) = 0, 
\end{align}
and, from \eqref{ecu13}, we define a new equation for $ \mathbf{Z}_{j}(\tau) $ by taking an arbitrary real number $ q_{j} $ such that:
\begin{align}\label{ecu34}
	\int_{0}^{2 \pi} \left\langle \mathbf{Z}_{0}(\tau), \mathbf{Z}_{j}(\tau) \right\rangle d \tau = q_{j}.
\end{align}

\begin{remark}\label{rmk1}
	The choice of the $ q_{j} $ is not critical for the final solution construction. Choosing different $ \{q_{j}\} $ sets will produce different $ \{\mathbf{Z}_{j} (\tau)\} $, $ \{\hat \lambda_{j}\} $ and $ \{\hat{T}_{j}\} $ solution sets, but all of them will eventually converge to the same solution with a different value of $ \varepsilon $. This is the reason for introducing the scaling factor $ \kappa(\varepsilon) $ in \eqref{ecu8} in the solution construction: different $ \{q_{j}\} $ sets produce different $ \kappa(\varepsilon) $, but the quantity $ \kappa (\varepsilon) \cdot \varepsilon^{2} $ in \eqref{ecu8} remains constant. Nevertheless,  it must also be noted that a sensible choice of the set $ \{q_{j}\} $ will have impact on the speed of the convergence of the series to the final solution. In our numerical experience, assuming the proximity of the final solution to the initial solution $ \mathbf{Z}_{0}(\tau) $, lower values of $ \{q_{j}\} $ produce a quicker convergence, so we use $ q_{j} = 0 \hspace{0.1 cm} \forall j $ in the included examples.
\end{remark}

Now we introduce the conditions for the system solvability.

\begin{dfn}\label{solvability} \textbf{System solvability }:
The system \eqref{ecu12}, \eqref{ecu13} and \eqref{ecu16} is formally solvable when the following conditions are met :	
	\begin{enumerate}[label=\roman*)]
		\item There exist unique $ \hat \lambda_{1} $ and $ \hat{T}_{1} $ such that:
		\begin{align}
			\nonumber \int_{0}^{2 \pi} \left\langle \hat{T}_{1} \mathbf{R}(\tau) + \hat \lambda_{1} \mathbf{S}(\tau) + M\mathbf{Z}_{0}(\tau) + N\mathbf{Z}_{0}(\tau - \hat \lambda_{0}) , \mathbf{w}_{i}(\tau) \right\rangle d \tau = 0 \hspace{0.5 cm} i=1, 2.
		\end{align}	
		
		\item It is verified that:
		\begin{align}\label{ecu36}
			\int_{0}^{2 \pi} \left\langle \hat{T} \mathbf{R}(\tau) + \hat \lambda \mathbf{S}(\tau), \mathbf{w}_{i}(\tau) \right\rangle d \tau = 0 \hspace{0.3 cm} i=1, 2 \Longrightarrow \hat \lambda = 0, \hat{T} = 0, 
		\end{align}		
	\end{enumerate}
where $ \{\mathbf{w}_{1}(\hat t), \mathbf{w}_{2}(\hat t)\} $ is the orthonormal basis of the null space of the adjoint operator $ N(L^{*}) $.
\end{dfn}

\begin{remark}\label{rmk2}
The two conditions given in Definition \ref{solvability} are sufficient for the system \eqref{ecu29}, \eqref{ecu33} and \eqref{ecu34} to have a unique solution $ \mathbf{Z}_{j}(\tau) $ for any $ q_{j} $ and $ \mathbf{P}_{j} $. Indeed, condition \eqref{ecu36} means that $ span(\{\mathbf{R}(\tau), \mathbf{S}(\tau) \}) \subset N(L^{*}) $, so, considering the orthogonality condition of $ R(L) $ with respect to the null space of the adjoint operator (i.e. $ R(L) = N(L^{*})^{\perp} $), we have:
\begin{itemize}
	\item For any $ \left\lbrace \hat{T}_{j}, \hat \lambda_{j} \right\rbrace \subset \mathbb{R} $ such that
	\begin{align}
		\nonumber \int_{0}^{2 \pi} \left\langle \hat{T}_{j} \mathbf{R}(\tau) + \hat \lambda_{j} \mathbf{S}(\tau), \mathbf{w}_{i}(\tau) \right\rangle d \tau = k_{i}, \hspace{0.1 cm} i=1, 2, 
	\end{align}
	there exists a unique $ \mathbf{Q}_{j}(\tau) \in N(L^{*}) $, $ \mathbf{Q}_{j}(\tau) = - \hat{T}_{j} \mathbf{R}(\tau) - \hat \lambda_{j} \mathbf{S}(\tau) $ such that 
	\begin{align}
	\nonumber \int_{0}^{2 \pi} \left\langle \mathbf{Q}_{j}, \mathbf{w}_{i}(\tau) \right\rangle d \tau = - k_{i}, \hspace{0.1 cm} i=1, 2.
	\end{align}
	Taking $ \mathbf{P}_{j}(\tau) = \mathbf{Q}_{j}(\tau) - M \mathbf{Z}_{j-1}(\tau) - N \mathbf{Z}_{j-1}(\tau-\hat \lambda_{0}) $ ensures the verification of the orthogonality condition.
	\item Conversely, for any $ \mathbf{P}_{j}(\tau) \in \mathbb{P} $ such that
	\begin{align}
	 \nonumber \int_{0}^{2 \pi} \left\langle \mathbf{P}_{j} + M \mathbf{Z}_{j-1}(\tau) + N \mathbf{Z}_{j-1}(\tau-\hat \lambda_{0}), \mathbf{w}_{i}(\tau) \right\rangle d \tau = k_{i}, \hspace{0.1 cm} i=1, 2, 
	\end{align}
	 there exist unique $ \left\lbrace \hat{T}_{j}, \hat \lambda_{j} \right\rbrace \subset \mathbb{R} $ such that 
	\begin{align}
	\nonumber \int_{0}^{2 \pi} \left\langle \hat{T}_{j} \mathbf{R}(\tau) + \hat \lambda_{j} \mathbf{S}(\tau), \mathbf{w}_{i}(\tau) \right\rangle d \tau = - k_{i}, \hspace{0.1 cm} i=1, 2, 
	\end{align}
	 verifying the orthogonality condition.
\end{itemize}
This, together with equations \eqref{ecu33} and \eqref{ecu34} (for any $ q_{j} $), ensure the uniqueness of $ \mathbf{Z}_{j}(\tau) $, which hence iteratively determine a unique solution for the system \eqref{ecu12}, \eqref{ecu13} and \eqref{ecu16}.
\end{remark}
	
To determine $ \mathbf{h}_j(\tau) $ in equation \eqref{ecu25}, we first must determine the unknown values for $ \hat \lambda_{j} $ and $ \hat{T}_{j} $. Using once again the orthogonality condition $ R(L) = N(L^{*})^{\perp} $, we get:
\begin{align}\label{ecu41} 
	\int_{0}^{2 \pi} \left\langle L\mathbf{Z}_{j} (\tau), \mathbf{w}_{i}(\tau) \right\rangle d \tau = \int_{0}^{2 \pi} \left\langle \mathbf{h}_{j} (\tau), \mathbf{w}_{i}(\tau) \right\rangle d \tau = 0 \hspace{10 mm} \text{for $ i=1, 2 $}.
\end{align}
The system \eqref{ecu41} generates two equations with the unknowns $ \hat \lambda_{j} $ and $ \hat{T}_{j} $, which provides the desired scalar values.


\subsubsection{Asymptotic validity of the solution}
The asymptotic validity of the solution calculated with the previous algorithm is ensured by the following theorem, whose proof can be found in the original reference.
\begin{thm}[Casal, Freedman \cite{Casal1980}]
	Under Assumption \ref{A1} and \ref{A2}, for an $ \varepsilon $ small enough, a sufficient condition for the system of equations \eqref{ecu12}, \eqref{ecu13} and \eqref{ecu16} to possess a unique and smooth solution $ \left( \hat \lambda(\varepsilon),  \hat{T}(\varepsilon), \mathbf{Z}(\tau, \varepsilon) \right) $, with $ \mathbf{Z}(\cdot, \varepsilon) $ periodic $ 2 \pi $, $ \hat \lambda(0) = \hat \lambda_{0} $, $ \hat{T}(0) = 2 \pi $ and $ \mathbf{Z}(\tau, 0) = \mathbf{z}_{0}(\tau) $, is that the system is \emph{formally solvable} in the sense of Definition \ref{solvability}.
	
	Under these conditions the solution $ \hat \lambda(\varepsilon),  \hat{T}(\varepsilon), \mathbf{Z}(\tau, \varepsilon) $ will possess asymptotic expansions of the form \eqref{ecu17}, \eqref{ecu18} and \eqref{ecu19} with the latter holding uniformly validly in $ \tau $ with $ \mathbf{Z}_{j}(\tau) $ $ 2 \pi-$periodic.
	
	In these expressions, the $ \left( \hat{\lambda}_{j}, \hat{T}_{j}, \mathbf{Z}_{j}(\tau) \right) $, are recursively the solutions of \eqref{ecu29}, \eqref{ecu33} and \eqref{ecu34}, whose existence is guaranteed by the conditions of \emph{formal solvability}.
\end{thm}

\section{Procedure to obtain a explicit solution of a desired order}\label{meth_ext}

\subsection{General expression for the solution of order $ j $}

Equation $ DDE_{Z_{j}} $ \eqref{ecu25} is the cornerstone to calculate the complete solution by an iterative process starting from $ Z_{0}(\tau) $ in \eqref{ecu17}. In the example included in \cite{Casal1980} the solution for these equations for orders 1 and 2 are directly presented, and there is not indication on how this is accomplished. Moreover, it is clear that for higher orders the $ DDE_{Z_{j}} $ can become highly complex. This problem is addressed in this section, where we provide a systematic approach to derive the solution up to the desired order of expansion for systems of any differential order, which can be efficiently implemented with standard computer tools.

For each order $ j $, after the calculation of $ \hat{T}_{j} $ and $ \hat \lambda_{j} $, and with $ \mathbf{Z}_{k}(\tau) $ known $ \forall k < j $, we look for solutions of \eqref{ecu29} that are periodic $ 2 \pi $. Equation \eqref{ecu25} generates an nonhomogeneous \textit{DDE} whose homogeneous version is coincident with \eqref{ecu6}, i.e:
\begin{align}\label{ecu45} 
	\nonumber L\mathbf{Z}_{j}(\tau) & = \mathbf{h}_{j}(\tau) = \mathbf{Z}_{j}'(t)- \mathbf{A}(\hat \lambda_{0}) \mathbf{Z}_{j}(\tau) - \mathbf{B}(\hat \lambda_{0}) \mathbf{Z}_{j}(\tau - \hat \lambda_{0}) \\
	& = \hat{T}_{j} \mathbf{R}(\tau) + \lambda_{j} \mathbf{S}(\tau) + M \mathbf{Z}_{j-1}(\tau) + N \mathbf{Z}_{j-1}(\tau - \hat \lambda_{0}) + \mathbf{P}_{j}(\tau), 
\end{align}
%
%
%
so the nonhomogeneous term $ \mathbf{h}_{j}(\tau) \in \mathbb{P} $, and is dependent on the following variables:
\begin{align}\label{ecu48} 
	\mathbf{h}_{j}(\tau) \equiv \mathbf{h}_{j}(\mathbf{Z}_{0}(\tau), ..., \mathbf{Z}_{j-1}(\tau), \hat \lambda_{0}, ..., \hat \lambda_{j}, \hat{T}_{0}, ..., \hat{T}_{j}).
\end{align}

To look for the periodic solutions $ \mathbf{Z}_{j}(\tau) $ we rely on the expression of the solution of a general linear \textit{DDE} initial value problem:
\begin{align}
	\nonumber \mathbf{u}(t)=\pmb{A}\mathbf{u}(t)+\pmb{B}\mathbf{u}(t-r)+\pmb{\Gamma}(t), \hspace{0.3 cm} t \geq 0, \hspace{0.3 cm} \mathbf{u}_{0}=\pmb{\phi}, 
\end{align}
where $ \pmb{A}, \pmb{B} \in \mathbb{R}^{n}\times \mathbb{R}^{n} $, $ \mathbf{u}_{0} $ is a continuous function $ \pmb{\phi}(t): [-r, 0) \longrightarrow \mathbb{R}^{n} $ representing the initial condition, and $ \pmb{\Gamma}(t): \mathbb{R}^{+} \rightarrow \mathbb{R}^n  $ is the nonhomogeneous term. The solution of this linear equation can be expressed in terms of the solution for the homogeneous problem and a particular solution for the nonhomogeneous problem with null initial condition (see \cite{Hale1977} and \cite{Smith2011}):
\begin{align}
	\nonumber \mathbf{u}(t)=\mathbf{u}(t, \pmb{\phi}, \pmb{\Gamma}) = \mathbf{u}(t, \pmb{\phi}, \mathbf{0}) + \mathbf{u}(t, \mathbf{0}, \pmb{\Gamma}).
\end{align}

Applying this principle to equation $ DDE_{Z_{j}} $ \eqref{ecu25}, we will calculate the periodic solutions of this equation as:
\begin{align}\label{ecu51} 
	\mathbf{Z}_{j}(\tau)= \overline {\mathbf{Z}}_{j}(\tau) + \widehat {\mathbf{Z}}_{j}(\tau),
\end{align}
where $ \overline {\mathbf{Z}}_{j}(\tau) $ is the general homogeneous solution and $ \widehat {\mathbf{Z}}_{j}(\tau) $ the particular solution for \eqref{ecu25}.

Regarding $ \overline {\mathbf{Z}}_{j}(\tau) $, it is easily calculated from the basis of $ 2 \pi $ periodic vectors $ \mathbf{v}_{i} $ for the null space of the operator $ L $ given in \eqref{ecu6}:
\begin{align}\label{ecu52} 
	& \overline {\mathbf{Z}}_{j}(\tau) = \sum_{i=1}^{2} \overline {c}_{ji}\mathbf{v}_{i}(\tau) =\overline{\mathbf{a}}_{j} \cos(\tau) + \overline{\mathbf{b}}_{j} \sin(\tau), 
\end{align}
where $ \overline c_{ji} \in \mathbb{R} $ and $ \overline {\mathbf{a}}_{j}, \overline {\mathbf{b}}_{j} \in \mathbb{R}^{n} $.

The calculation of the particular solution $ \widehat {\mathbf{Z}}_{j}(\tau) $ is more involved. We start by noting the structure of the nonhomogeneous term $ \mathbf{h}_{j}(\tau) $ in \eqref{ecu48} by the following proposition.
\begin{prop}\label{Prop1}
\item The function $ \mathbf{h}_{j}(\tau) $ can be expressed as a truncated Fourier series up to the index $ j+1 $, that is:
\begin{align}\label{ecu53} 
\mathbf{h}_{j}(\tau) = \pmb{\alpha}_{j0}+\sum_{k=1}^{j+1} (\pmb{\alpha}_{jk}\cos(k \tau) + \pmb{\beta}_{jk}\sin(k \tau)), 
\end{align}
where $ \pmb{\alpha}_{jk}, \pmb{\beta}_{jk} \in \mathbb{R}^n $.
\end{prop}

\begin{proof}
A solution $ \mathbf{Z}_{0}(\tau) $ to the homogeneous equation:	
\begin{align}
	\nonumber \mathbf{Z}_{0}'(\tau)=\mathbf{A}(\hat \lambda_{0}) \mathbf{Z}_{0}(\tau) + \mathbf{B}(\hat \lambda_{0}) \mathbf{Z}_{0}(\tau-\hat \lambda_{0})
\end{align}
can be built from the base $ \mathbf{v}_{i}(\tau) $:
\begin{align}\label{ecu55} 
	\mathbf{Z}_{0}(\tau) = \sum_{i=1}^{2} C_{0i}\mathbf{v}_{i}(\tau), 
\end{align}
with $ C_{0i} \in \mathbb{R} $.

From Assumption \ref{A1}, the linear equation \eqref{ecu2} presents a pair of pure imaginary roots $ \pm j \omega_{0} $, so a basis $ \{ \mathbf{\widetilde{v}}_{i}(t) \}, i=1..2 $ for the periodic solutions of equation \eqref{ecu2} can have the form:
\begin{align}
	\nonumber \mathbf{\widetilde{v}}_{i}(t) = \pmb{\alpha}_{i} \cos(\omega_{0}t) + \pmb{\beta}_{i} \sin(\omega_{0}t), 
\end{align}
where $ \pmb{\alpha}_{i} , \pmb{\beta}_{i} \in \mathbb{R}^n $ can be chosen so that the basis $ \{ \mathbf{\widetilde{v}}_{i}(t) \} $ is orthonormal.
Then an orthonormal basis for the solutions of equation \eqref{ecu4} can have the form:
\begin{align}\label{ecu57} 
	\mathbf{v}_{i}(\hat t) = \mathbf{\widetilde{v}}_{i}\left( \frac{\hat{t}}{\omega_{0}}\right) = \pmb{\alpha}_{i} \cos(\hat t) + \pmb{\beta}_{i} \sin(\hat t).
\end{align}
So, \eqref{ecu55} and \eqref{ecu57} yield that $ \mathbf{Z}_{0} $ is a first order trigonometric polynomial, that is:
\begin{align}\label{ecu58} 
	\mathbf{Z}_{0}(\tau) = \sum_{i=1}^{2} C_{0i}\mathbf{v}_{i}(\tau)= \sum_{i=1}^{2} C_{0i}(\pmb{\alpha}_{i} \cos(\tau) + \pmb{\beta}_{i} \sin(\tau)).
\end{align}
Recalling \eqref{ecu45}, the function $ \mathbf{h}_{j}(\tau) $ has the following form in each step:
\begin{align}\label{ecu59} 
	\mathbf{h}_{j}(\tau)=\hat{T}_{j} \mathbf{R}(\tau) + \lambda_{j} \mathbf{S}(\tau)
	+ M \mathbf{Z}_{j-1}(\tau) + N \mathbf{Z}_{j-1}(\tau - \hat \lambda_{0}) + \mathbf{P}_{j}(\tau).
\end{align}
We note, for each component of $ \mathbf{h}_{j}(\tau) $ in \eqref{ecu59}, the following properties:
\begin{enumerate}[label=(\alph*)]
	\item The two first terms are a linear operator on $ \mathbf{Z}_{0}(\tau) $, $ \mathbf{Z}_{0}(\tau - \hat \lambda_{0}) $ and their first derivatives (equations \eqref{ecu30} and \eqref{ecu31}), thus yielding a first order trigonometric polynomial.
	\item The third and fourth terms of \eqref{ecu59} increase the trigonometric order of $ \mathbf{Z}_{j-1}(\tau) $ in 1  by means of the bilinear forms $ M $ \eqref{ecu27} and $ N $ \eqref{ecu28} applied to $ \mathbf{Z}_{j-1}(\tau) $ and $ \mathbf{Z}_{j-1}(\tau - \hat \lambda_{0}) $ respectively.
	\item $ \mathbf{P}_{j}(\tau) $ is a trigonometric polynomial of order $ j + 1 $ by virtue of cross products of $ \mathbf{Z}_{k}, k < j $ and their derivatives, plus the expansion term of order $ j $ of \eqref{ecu23}.
\end{enumerate}

Finally, taking into account that the order 0 solution $ \mathbf{Z}_{0}(\tau) $ has by \eqref{ecu58} a trigonometric order $ 1 $, the preceding properties $ ((a) - (c)) $ ensure that equation \eqref{ecu53} is verified for all orders of $ j $, thus proving the proposition. 
\end{proof}

The next Corollary enable us to calculate $ \widehat {\mathbf{Z}}_{j}(\tau) $:

\begin{cor}\label{Cor1}
	The particular solution of \eqref{ecu25} $ \widehat {\mathbf{Z}}_{j}(\tau) $ is a trigonometric polynomial of degree $ j+1 $, which can be expressed as:
	\begin{align}\label{ecu60}
	\widehat{\mathbf{Z}}_{j} (\tau)= 	\widehat{\mathbf{a}}_{j0} + \sum_{k=1}^{j+1} (\widehat{\mathbf{a}}_{jk} \cos(k \tau) + \widehat{\mathbf{b}}_{jk} \sin(k \tau)), 
	\end{align}
	where $ \widehat {\mathbf{a}}_{jk}, \widehat {\mathbf{b}}_{jk} \in \mathbb{R}^n $.
	
	\noindent Moreover:
	\begin{align}\label{ecu61}
		\widehat{\mathbf{Z}}_{j}' (\tau)= \sum_{k=1}^{j+1} (-\widehat{\mathbf{a}}_{jk} k \sin(k \tau) + \widehat{\mathbf{b}}_{jk} k \cos(k \tau)).
	\end{align}
	Then the particular solution can be calculated by solving the system obtained by introducing the expressions \eqref{ecu60} for $ \mathbf{Z}_{j} (\tau) $ and \eqref{ecu61} for $ \mathbf{Z}_{j}' (\tau) $ in the identity $ DDE_{Z_{j}} $ \eqref{ecu25}. Equating the coefficients of the same trigonometric term in each of the $ n $ components produces a linear system of order $ n(2j+3) $ which provides the desired solution for the coefficients $ \widehat{\mathbf{a}}_{jk} $ and $ \widehat{\mathbf{b}}_{jk} $.
\end{cor}
\begin{proof}
	$ \widehat {\mathbf{Z}}_{j}(\tau) $ and $ \widehat {\mathbf{Z}}'_{j}(\tau) $ must have the same trigonometric degree, so by \eqref{ecu25} and Proposition \ref{Prop1} that degree must be $ j+1 $. The rest is a consequence of this fact.
\end{proof}

%
%
%

\subsection{Solution $ \mathbf{Z}_{j} (\tau) $ of $ DDE_{Z_{j}} $ \eqref{ecu25}}\label{complete_sol_i}

To determine $ \mathbf{Z}_{j} (\tau) $ by \eqref{ecu51} we need to have $ \overline {\mathbf{Z}}_{j}(\tau) $ and $ \widehat {\mathbf{Z}}_{j}(\tau) $. The particular solution $ \widehat {\mathbf{Z}}_{j}(\tau) $ has been already obtained by Corollary \ref{Cor1}, so to uniquely determine the solution to the homogeneous equation $ \overline {\mathbf{Z}}_{j}(\tau) $,  we get the coefficients $ \overline c_{ji} $ in \eqref{ecu52} for every iteration by applying the supplementary conditions for orthogonality and phase with respect to the order 0 solution $ \mathbf{Z}_{0} $.

\noindent To get $ \mathbf{Z}_j^1(0) = 0 $ (see \eqref{ecu33}), we use the expressions \eqref{ecu52} and \eqref{ecu60}, such that:
\begin{align}\label{ecu66}
	0 = \mathbf{Z}_{j}^{1}(0) = \widehat{\mathbf{Z}}_{j}^{1} (0) + \overline{\mathbf{Z}}_{j}^{1}(0)= \widehat{\mathbf{a}}_{j0}^{1} + \sum_{k=1}^{j+1} \widehat{\mathbf{a}}_{jk}^{1} + \overline{c}_{j1} \mathbf{v}_{1}^{1}(0) + \overline{c}_{j2} \mathbf{v}_{2}^{1}(0).
\end{align}
From \eqref{ecu34}, we choose $ q_{j}=0 $ for $ j \geq 1 $:
\begin{align}\label{ecu67} 
	0 = \int_{0}^{2 \pi} \left\langle \mathbf{Z}_{j} (\tau), \mathbf{Z}_{0}(\tau) \right\rangle d \tau = 	\int_{0}^{2 \pi} \left\langle \widehat{\mathbf{Z}}_{j}(\tau) + \overline{\mathbf{Z}}_{j}(\tau), \mathbf{Z}_{0}(\tau) \right\rangle d \tau.
\end{align}
The system \eqref{ecu66} - \eqref{ecu67} provides the solution for $ \overline{c}_{ji} , i=1, 2 $ that uniquely determine $ \overline{\mathbf{Z}}_{j}(\tau) $. This finishes the calculation of $ \mathbf{Z}_{j} (\tau) $.


\subsection{Complete solution}
  Following this procedure, for the final solution $ \mathbf{Z}(\tau, \varepsilon) $ two series expansions are carried out simultaneously, one asymptotic in powers of $ \varepsilon $, and the other Fourier-type in periods multiple of $ 2 \pi $.
\begin{align}\label{ecu67b}
	&\mathbf{Z}(\tau, \varepsilon) = \sum_{j=0}^{N} \mathbf{Z}_{j}(\tau)\varepsilon^{j}, \\ 
	&\mathbf{Z}_{j} (\tau)= 	\mathbf{a}_{j0} + \sum_{k=1}^{j+1} (\mathbf{a}_{jk} \cos(k \tau) + \mathbf{b}_{jk} \sin(k \tau)), 
\end{align}
that, together with \eqref{ecu51}, \eqref{ecu52} and \eqref{ecu60} yield:
\begin{align}
		\nonumber \mathbf{a}_{j0} &= \widehat {\mathbf{a}}_{j0}, & & & \text{for }k&=0,  \\
	  \nonumber \mathbf{a}_{j1} &= \widehat {\mathbf{a}}_{j1} + \overline {\mathbf{a}}_{j}, & \mathbf{b}_{j1} & = \widehat {\mathbf{b}}_{j1} + \overline {\mathbf{b}}_{j}, & \text{for }k&=1, \\
		\nonumber \mathbf{a}_{jk} &= \widehat {\mathbf{a}}_{jk}, & \mathbf{b}_{jk} & = \widehat {\mathbf{b}}_{jk}, & \text{for }k&>1.
\end{align}

At last, the explicit approximate solution to the original equation up to order $ N $ is recovered as:
\begin{align}
	\nonumber \mathbf{x}_{N}(t) = \tilde \varepsilon \left( \mathbf{Z}_{0}\left( \dfrac{2 \pi}{\hat{T}(\tilde \varepsilon)} \omega_{0} t\right) + \sum_{j=1}^{N} \mathbf{Z}_{j}\left( \dfrac{2 \pi}{\hat{T}(\tilde \varepsilon)} \omega_{0} t\right) \tilde \varepsilon^{j} \right),	
\end{align}
where $ \tilde \varepsilon $ is the solution to the equation:
\begin{align}\label{ecu71}	
	\lambda = \frac{1}{\omega_{0}} \sum_{j=0}^{N} \hat \lambda_{j}\tilde \varepsilon^{j}.	
\end{align}
\subsection{Error estimation}
Strict error bounds for approximate solutions in systems of differential equations are difficult to calculate. For ODEs and PDEs the bounds are usually based on specific characteristics of the differential equation, and no general methodology exist. If we denote the differential problem as $ \mathbf{x}' = \mathbf{f}(t, \mathbf{x}) $ and $ \mathbf{x}_{N} $ as the approximation of order $ N $, it is usually searched a relationship of the error function calculated in some norm $ \|\mathbf{x} - \mathbf{x}_{N} \| $, with the \textit{residual} function $ \mathbf{R}_{N} $ that measures how well $ \mathbf{x}_{N} $ satisfies the differential equation:
\begin{align}
	\nonumber \mathbf{R}_{N}(t) = \mathbf{x}_{N}'(t) - \mathbf{f}(t, \mathbf{x}_{N}(t)).
\end{align}
If the problem is well-posed, then it is considered that the residual represents a good measure of the deviation of the approximation from the real solution (see, for instance, \cite{Ferguson1972}).

In the specific case of RFDEs of discrete type, the difficulty of the problem is compounded by the delay terms. There exist recent works that allow to calculate an error bound for approximate periodic solutions in differential delay equations by the estimation of a set of parameters and the verification of three inequalities \cite{Gilsinn2009}. Nevertheless, these calculations require a significant amount of computing effort, that may be beyond the effort to calculate the approximation itself. Therefore, given this difficulty and considering the regularity of the problem we are addressing, we use in this paper for the measure of the relative error the supremum norm of the residual in $ [0, \hat{T}_{N}] $, where $ \hat{T}_{N} $ is the period of the approximate solution, i.e.:
\begin{align}\label{ecu73} 
	r_{r} = \frac{\|\mathbf{R}_{N}(t)\|_{\infty}}{\|\mathbf{x}_{N}'(t)\|_{\infty}} = \frac{\sup\left\lbrace |\mathbf{x}_{N}'(t) - \mathbf{f}(t, \mathbf{x}_{N}(t), \mathbf{x}_{N}(t - \lambda))| : t \in [0, \hat{T}_{N}] \right\rbrace}{\sup\left\lbrace |\mathbf{x}_{N}'(t)| : t \in [0, \hat{T}_{N}] \right\rbrace} .
\end{align}

In addition to this measure of residual error, in the  examples in next section we will also compare the approximate solutions with numerical integrations obtained with the solver routine for \textit{DDEs} of \textsf{Maple 2023} (\textsf{dsolve}). This numerical solution, denoted $ \tilde {\mathbf{x}}(t) $, becomes quasi-periodic with a period $ \tilde T $ after a long enough initial time $ t_{0} $. In addition to numerically verify quasi-periodicity, the precise value of $ t_{0} $ is further determined by the conditions $ \tilde{\mathbf{x}}^{1}(t_{0})=\mathbf{x}^{1}_{N}(0)=0  $ and $ \sgn\left( d\tilde{\mathbf{x}}^{1}/dt (t_{0})\right) = \sgn\left( d\mathbf{x}^{1}_{N}/dt (0)\right)  $. Our estimation for the relative error of the calculated solution $ \mathbf{x}_{N}(t) $ with respect to the numerical solution $ \tilde{\mathbf{x}}(t) $ is thus based on this time reference:
\begin{align}\label{ecu74} 
	e_{r} = \frac{\|\tilde {\mathbf{x}}(t_{0}+t) - \mathbf{x}_{N}(t)\|_{\infty}}{\|\tilde {\mathbf{x}}(t_{0}+t)\|_{\infty}} = \dfrac{\sup\left\lbrace |\tilde {\mathbf{x}}(t_{0}+t) - \mathbf{x}_{N}(t)| : t \in [0, \tilde T] \right\rbrace}{\sup\left\lbrace |\tilde {\mathbf{x}}(t_{0}+t)| : t \in [0, \tilde T] \right\rbrace}.
\end{align}
The advantage of this approach is that in the case that the numerical method is error controlled, in the sense that $ \|{\mathbf{x}}(t_{0}+t) - \tilde {\mathbf{x}}(t_{0}+t)\|_{\infty} \leq \delta $, it provides a hard bound for the absolute error of the calculated approximation $ \mathbf{x}_{N}(t) $, that is, $ \|{\mathbf{x}}(t_{0}+t) - \mathbf{x}_{N}(t)\|_{\infty} \leq \|{\mathbf{x}}(t_{0}+t) - \tilde {\mathbf{x}}(t_{0}+t)\|_{\infty} + \|\tilde {\mathbf{x}}(t_{0}+t) - \mathbf{x}_{N}(t)\|_{\infty} = \delta +   \|\tilde {\mathbf{x}}(t_{0}+t) - \mathbf{x}_{N}(t)\|_{\infty} $. In this way we can determine with the aid of the numerical solution $ \tilde {\mathbf{x}}(t) $ the minimum expansion order $ N $ that provides a maximum error in a specific scenario.

\section{Numerical examples}

\subsection{A car following model for road traffic analysis}
Many road traffic studies are focused on the determination of the conditions that generate congestion and traffic jams, for which two different approaches are usually used to perform these analyses: the \textit{microscopic model}, where traffic is seen as individual interacting particles, and the \textit{macroscopic model}, where traffic is seen as a compressible fluid. 

A typical micro model is the \textit{follow-the-leader} or \textit{car-following} model, where the driver's actions in a vehicle $ n $ are induced by the relative motion of the vehicle $ n-1 $ ahead. Most important among these models are those that consider the influence of the inter-vehicle distance ($ X_{n-1}(t) - X_{n}(t) $) and the difference between their respective velocities ($ X'_{n-1}(t) - X'_{n}(t) $) \cite{Gazis1959, Gazis1961, Edie1961, Newell1961, Bando1995}. We will focus on a new car-following model (\textit{\textbf{nDDE}}) recently introduced \cite{Padial2022} that present advantages from the point of view of the phenomenological explanation, and for which, conveniently, a complete analysis of the Hopf bifurcation that it presents is also provided.

This new model represents the situation of two vehicles in a single lane (denoted as $ 0 $ for the leading car and $ 1 $ for the following car). It considers that the leading car has a constant velocity $ X'_{0}(t)=v_{0} $, and the following variables for the distance between the cars ($ s(t) $) and the relative velocity ($ s'(t) $):
\begin{align}
	\nonumber &s(t)=(X_{0}(t)-X_{1}(t)), &s'(t)=(v_{0}-X_{1}'(t)).	
\end{align}
We then have $ s''(t)=-X''_{1}(t) $. The authors in \cite{Padial2022} express the retarded dependence of the acceleration as a function of $ s(t) $ and $ s'(t) $. In this way we have the system (\textbf{\textit{nDDE}}):
\begin{align}\label{ecu82}
	-s''(t + \lambda)=X''_{1}(t+\lambda)=g(s(t), s'(t)), 
\end{align}
where $ g $ is the following sigmoid function:

\begin{align}\label{ecu83}
	g(s, s')=a-\frac{(a+b)}{1+\frac{b}{a}e^{d(s-M+Ks')}}, \space \forall(s, s') \in \mathbb{R}^{2}.
\end{align}
In \eqref{ecu82} and \eqref{ecu83} $ \lambda $ is the delay or reaction time, $ M $ is the equilibrium distance, $ d $ is a parameter to model the intensity of the response of the car-driver ensemble, and $ K $ is a parameter to model the driver's response according to the safe distance and the perceived relative velocity.

Reference \cite{Padial2022} provides an analysis of the dynamical characteristics of the Hopf bifurcation, showing the change from stability to oscillatory solutions in equilibrium depending on the value of the delay term. It also characterizes the \textit{nDDE} equation as a Retarded Functional Differential Equation (RFDE) \cite{Hale1977}, and studies the structure of the solutions produced by this bifurcation and its impact on the traffic process. As an application of the methodology introduced in this paper, next we use it to calculate the explicit expression for the periodic solutions arising near the bifurcation point for this model.
\subsubsection{System of equations and equilibrium solution}
We start by redefining the variables so that the equilibrium point is $ \mathbf{0} $, defining $ x_{1} = s - M $ and $ x_{2} = s' $, and transforming equation \eqref{ecu82} into a system of first order \textit{DDEs}:
\begin{align}\label{ecu84}
	\nonumber & x_1'(t) = x_{2}(t), & \\
	& x_{2}'(t) = -a + \frac{(a+b)}{1+\frac{b}{a}e^{d(x_1(t-\lambda)+Kx_{2}(t - \lambda))}}.
\end{align}
Next, we consider the linearized system:
\begin{align}\label{ecu85} 
	\nonumber & x_1'(t) = x_{2}(t), & \\
	& x_{2}'(t) = -D * \left( x_1(t-\lambda)+K x_{2}(t - \lambda) \right), 
\end{align}
where $ D = d\frac{ab}{a+b} $.
This system is characterized by the following matrices (as in equation \eqref{ecu2}):
\begin{align}
	\nonumber \mathbf{P} = \left(\begin{matrix}
		0 & 1\\
		0 & 0
	\end{matrix}\right), \hspace{2 cm} \mathbf{Q} = \left(\begin{matrix}
		0 & 0\\
		-D & -DK
	\end{matrix}\right).
\end{align}

In order to calculate the critical delay $ \lambda_{0} $ that produces the periodic solution of period $ \omega_{0} $ corresponding to the Hopf bifurcation, we proceed in the usual way to derive the characteristic equation by trying an exponential function as the eigenfunction of the system. Introducing a general solution of the form $ \pmb{\alpha} e^{j \omega_{0}t} $ , with $ \pmb{\alpha} \in \mathbb{R}^{2} $ constant, we get the characteristic equation:
\begin{align}
	\nonumber \omega_{0}^2 = D ( 1 + j K \omega_{0}) e^{j \omega_{0} \lambda_{0}}, 
\end{align}
which, in turn, generates the following equations corresponding to the real and imaginary parts:
\begin{flalign}\label{ecu97}
	K = \frac{\tan(\omega_{0} \lambda_{0})}{\omega_{0}}, \hspace{2 cm}
	D = \omega_{0}^2 \cos(\omega_{0} \lambda_{0}).
\end{flalign}
These two transcendental equations define the bifurcation point $ (\omega_{0}, \lambda_{0}) $.

\paragraph{Solution bases}
We get the following solution basis in terms of the scaled time $ \hat t = \omega_{0} t $:
\begin{align}
	&\nonumber \mathbf{v}_{1}(\hat t) = \frac{1}{\sqrt{\pi(1 + \omega_{0}^2)}} (\cos(\hat t), -\omega_{0} \sin(\hat t)),
	&\nonumber \mathbf{v}_{2}(\hat t) = \frac{1}{\sqrt{\pi(1 + \omega_{0}^2)}} (\sin(\hat t), \omega_{0} \cos(\hat t)).
\end{align}
For the adjoint system, we get the following basis in $ \hat t $:
\begin{flalign}
	\nonumber \mathbf{w}_{1}(\hat t) = \frac{\cos(\omega_{0} \lambda_{0})}{\sqrt{\pi(1 + \omega_{0}^2 \cos^2(\omega_{0} \lambda_{0}))}} (\omega_{0} \cos(\hat t), \tan(\omega_{0} \lambda_{0})\cos(\hat t) - \sin(\hat t)), \\
	\nonumber \mathbf{w}_{2}(\hat t) = \frac{\cos(\omega_{0} \lambda_{0}) }{\sqrt{\pi(1 + \omega_{0}^2 \cos^2(\omega_{0} \lambda_{0}))}} (\omega_{0} \sin(\hat t), \cos(\hat t) + \tan(\omega_{0} \lambda_{0})\sin(\hat t)).
\end{flalign}

\subsubsection{Explicit expression of periodic solutions}
From now on we will focus on a particular case of system \eqref{ecu84}, using the same physical parameters as in \cite{Padial2022}, summarized in Table \ref{Tab-2}.

\begin{table}[!htp]
	\centering
	\begin{tabular}{clc}
		\toprule
		Param & Description & Value \\
		\midrule
		$ a $ & Maximum acceleration in $ nDDE$& $ 2.0576 \text{  } m/s^{2} $ \\
		$ b $ & Maximum deceleration & $ 1.5677 \text{  } m/s^{2} $ \\
		$ v_{0} $ & Leader velocity & $ 22.2222 \text{  } m/s $ \\
		$ M $ & Safe distance for $ v_{0} $ & $ 44.4444 \text{  } m $ \\
		$ d $ & Car + driver response intensity & $ 0.1124 $ \\
		$ K $ & Driver sensitivity to $ s' $ & $ 11.3890 \text{  } s $ \\
		\bottomrule
	\end{tabular}
	\caption{\textit{nDDE} parameters for the evaluated scenario \cite{Padial2022}.}\label{Tab-2} 
\end{table}
Under these conditions, the bifurcation point can be calculated numerically from \eqref{ecu97}:
\begin{align}
 \nonumber \omega_{0} = 1.1424, \hspace{2 cm} \lambda_{0} = 1.3079, 
\end{align}
obtaining the following orthonormal base for the solutions of the linearized system \eqref{ecu85}:
\begin{align}
	&\nonumber \mathbf{v}_{1}(\hat t) = (0.3716 \cos(\hat t), -0.4245 \sin(\hat t)), 
	&\nonumber \mathbf{v}_{2}(\hat t) = (0.3716 \sin(\hat t), 0.4245 \cos(\hat t)), 
\end{align}
and the base for the adjoint system:
\begin{align}
	&\nonumber \mathbf{w}_{1}(\hat t) = (0.0492 \cos(\hat t), -0.0431 \sin(\hat t) + 0.5604 \cos(\hat t)), \\
	&\nonumber \mathbf{w}_{2}(\hat t) = (-0.0492 \sin(\hat t), -0.0431 \cos(\hat t) - 0.5604 \sin(\hat t)).
\end{align}

We then use these values to derive the approximations to the solutions of $ \mathbf{Z}(\tau, \varepsilon), \hat \lambda(\varepsilon) $ and $ \hat{T}(\varepsilon) $ by iteratively calculating the terms $ \mathbf{Z}_{j} (\tau) $, $ \hat \lambda_{j}$ and $ \hat{T}_{j} $ as described in Sections \ref{PL_method} and \ref{meth_ext}. We have performed the automation of the symbolic and numerical calculations by means of a program developed on the \textsf{Maple 2023} platform.

First, using the basis of the linear system $ \{\mathbf{v}_{1}, \mathbf{v}_{2}\} $, we start by selecting the fundamental solution $ \mathbf{Z}_{0}(\tau) $ so that it is compatible with the phase condition \eqref{ecu22}. This leads us to the choice $ \mathbf{Z}_{0}(\tau) = 	\mathbf{v}_{2}(\tau) $. Using this as our first solution, we iteratively calculate the subsequent orders following the described methodology. The results are presented (up to order 5) in different tables. 
\begin{table}[ht]
	\centering
	\begin{tabular}{ccccccc}
		\hline
		Coef. & o0 & o1 & o2 & o3 & o4 & o5 \\ 
		\hline
		$ \hat \lambda_{j} $ & 1.4940 &   0.0000 & 0.1666 &   0.0000 & 0.0387 & 0.0013 \\ 
		$\hat{T}_{j} $ & $ 2 \pi $ &   0.0000 & 0.7465 &   0.0000 & 0.1814 & 0.0056 \\ 
		\hline
	\end{tabular}
	\caption{Coefficients for series expansion of $ \hat \lambda $ and $ \hat{T} $.} 
	\label{coef_lt}
\end{table}
\begin{table}[ht]
	\centering
	\addtolength{\tabcolsep}{-0.4em}
	\begin{tabular}{|c|cccccc|cccccc|}
		\multicolumn{1}{c}{} &  \multicolumn{6}{c}{\large $ x_{1} $ } & \multicolumn{6}{c}{\large $ x_{2} $ } \\ 
		\hline
		Trig. order & o0 & o1 & o2 & o3 & o4 & o5 & o0 & o1 & o2 & o3 & o4 & o5 \\ 
		\hline
		Const & 0.0000 & -3.5250 & 0.0000 & -6.4153 & -0.4194 & -13.7918 & 0.0000 & 0.0000 & 0.0000 & 0.0000 & 0.0000 & 0.0000 \\ 
		$ \cos(\tau) $ & 2.3349 & 0.0000 & 0.9868 & 0.0000 & 1.1873 & 0.1174 & 0.0000 & -4.0910 & -0.1035 & -3.9406 & -0.5446 & -7.6870 \\ 
		$ \sin(\tau) $ & 0.0000 & 3.5811 & 0.0906 & 6.1228 & 0.5443 & 13.3860 & 2.6673 & 0.0000 & -0.8638 & 0.0000 & -1.0393 & -0.1028 \\ 
		$ \cos(2\tau) $ & & -0.0295 & 0.1722 & 0.0741 & 0.3616 & 0.2281 & & 0.1282 & 0.2068 & -0.2905 & 0.0943 & -0.8224 \\ 
		$ \sin(2\tau) $ & & -0.0561 & -0.0905 & 0.0852 & -0.1088 & 0.3909 & & -0.0674 & 0.3933 & 0.2196 & 0.5325 & 0.4340 \\ 
		$ \cos(3\tau) $ & & & 0.0461 & 0.0003 & -0.3647 & -0.0313 & & & 0.0002 & -0.7267 & -0.0311 & 0.8003 \\ 
		$ \sin(3\tau) $ & & & -0.0001 & 0.2120 & 0.0090 & -0.0752 & & & 0.1579 & 0.0009 & -1.3679 & -0.1078 \\ 
		$ \cos(4\tau) $ & & & & -0.0041 & 0.0293 & 0.0632 & & & & 0.0218 & 0.1153 & -0.3241 \\ 
		$ \sin(4\tau) $ & & & & -0.0048 & -0.0252 & 0.0674 & & & & -0.0188 & 0.1339 & 0.3028 \\ 
		$ \cos(5\tau) $ & & & & & 0.0030 & -0.0007 & & & & & -0.0005 & -0.1333 \\ 
		$ \sin(5\tau) $ & & & & & 0.0001 & 0.0233 & & & & & 0.0174 & -0.0042 \\
		$ \cos(6\tau) $ & & & & & & -0.0006 & & & & & & 0.0038 \\ 
		$ \sin(6\tau) $ & & & & & & -0.0006 & & & & & & -0.0040 \\ 
		\hline
	\end{tabular}
	\caption{Coefficients for space $ x_{1} $ and relative velocity  $ x_{2} $.} 
	\label{coef}
\end{table}
Table \ref{coef_lt} contains the coefficients for the series expansion of $ \hat \lambda $ and $ \hat{T} $, up to order 5. Table \ref{coef} contain the coefficients for each trigonometric order for the space $ x_{1}(t) $ and the relative velocity $ x_{2}(t) $ series expansions. Notice that Table \ref{coef_lt} for the delay and period is provided in terms of the scaled time $ \hat t = \omega_{0} t $, while Table \ref{coef} for the periodic solutions are provided in terms of the scaled time $ \tau = (2 \pi/\hat{T}(\varepsilon)) \cdot \hat t = (2 \pi/\hat{T}(\varepsilon)) \cdot \omega_{0} t $. It is worth noting that even when we get null values for the parameters $ \hat \lambda_{i} $ and $ \hat{T}_{i} $, as is the case for order 1 and 3, there are still nonzero terms for those orders in the series expansion of the solution.

At this point we are already able to generate valuable insights based on these coefficients. We can, for instance, draw the bifurcation diagram for the different components of the solution by representing the maximum and minimum of the oscillations for each value of the delay. This is shown in Figure \ref{Gr0_nDDE} below. Notice that this diagram is generated almost instantly from the series expansion obtained, with no need to perform numerical integrations for each value of $ \lambda $, as is customary when no explicit expression for the solution exists.

\begin{figure}[!h]\centering
	\includegraphics[width=0.6 \textwidth]{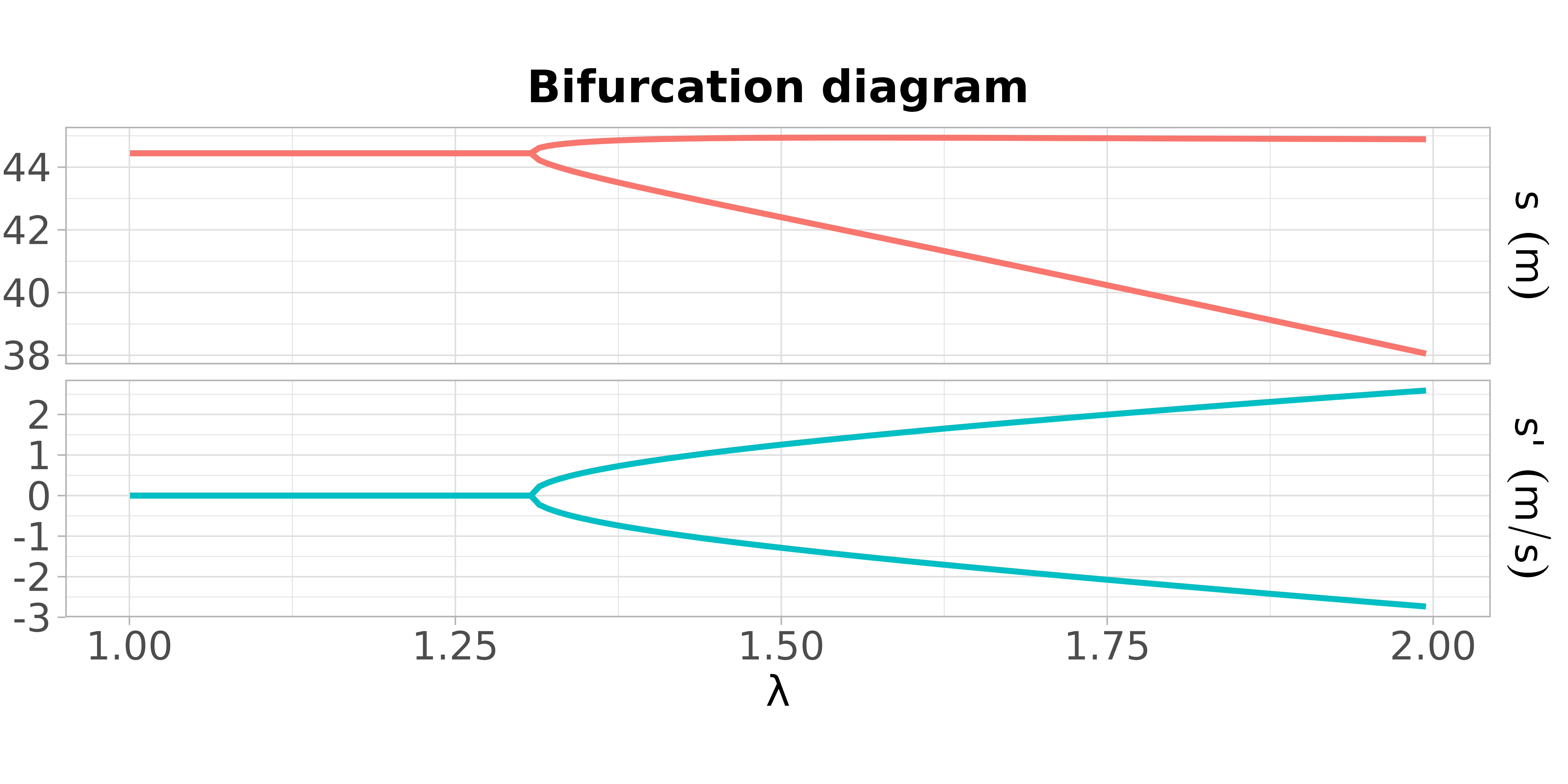}
	\caption{Bifurcation diagram of $ s $ and $ s' $ in the \textit{nDDE} model \eqref{ecu83}.}
	\label{Gr0_nDDE}
\end{figure}

More detailed results for specific values of delay follow. Given that the periodic solution for values of delay that are very close to the bifurcation value $ \lambda_{0} = 1.3079 $ have very small amplitude, we test delay values sensibly greater to get oscillations with significant amplitudes. The results compare the approximation obtained with the Poincaré-Lindstedt method up to order 8 (referred to as \textit{P-L} or $ \mathbf{x}_{8} $ in the graphs) to the solution calculated by the numerical integration (referred to as  \textit{Num}). The Maple routine used is \textsf{dsolve}. In this case, using the algorithm \textsf{rkf45} with parameters \textsf{abserr} = \textsf{relerr} = $ 10^{-9}$ (for details of this solver, see \cite{maple_dde}). The numerical integration results are calculated with an initial condition $ x_{1}(0) = 20 \hspace{0.2 cm} m $ and $ x_{2}(0) = 17.22 \hspace{0.2 cm} m/s $, taking the steady state after a sufficiently large initial time $ t_{0} > 1000 s $, while the explicit calculation by Poincaré-Lindstedt method is directly calculated up to the order 8. 

\begin{figure}[h!]\centering
	\begin{subfigure}{0.49\textwidth}
		\includegraphics[width=\textwidth]{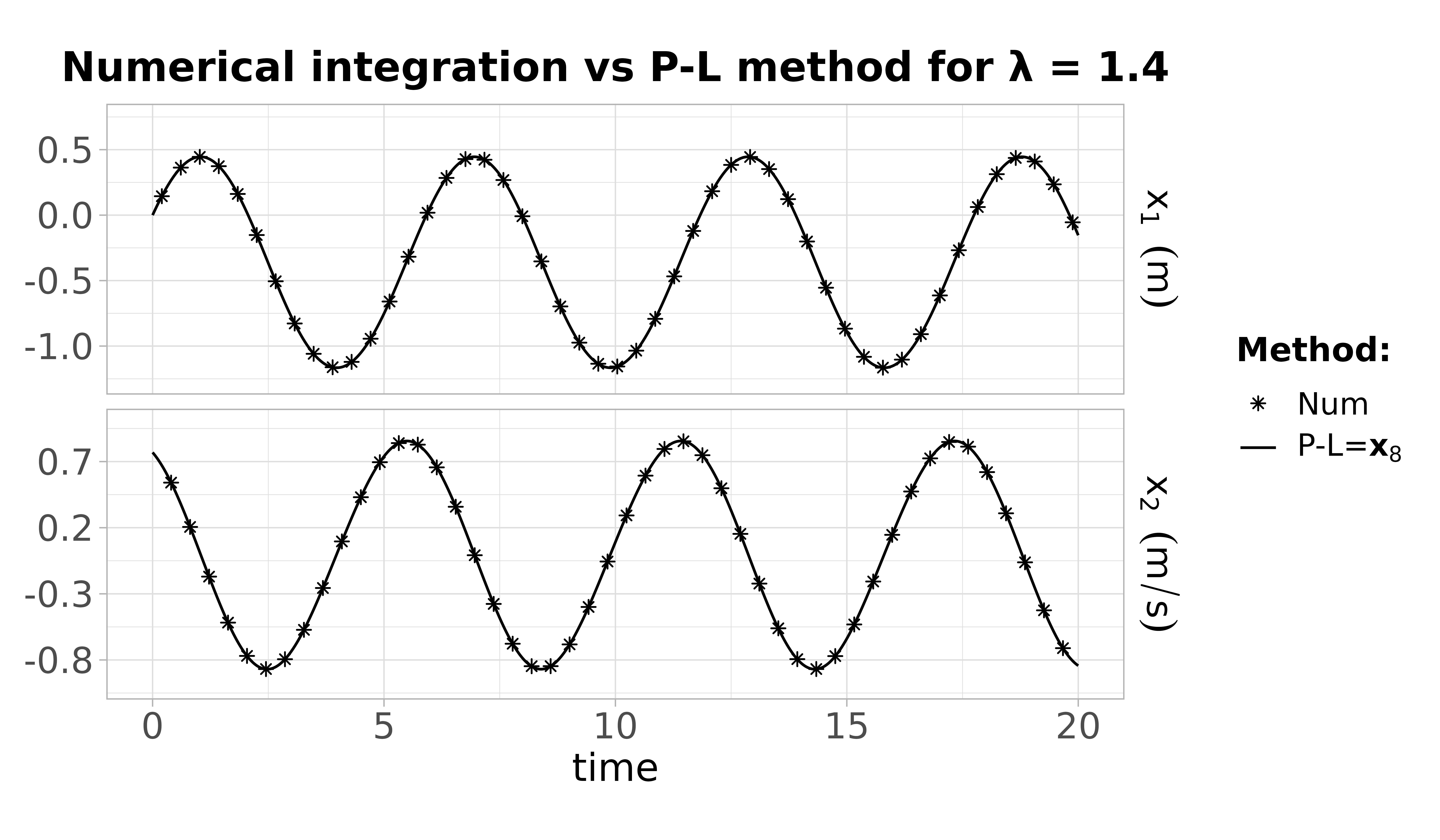}
		\caption{Numerical solution $ \tilde {\mathbf{x}}(t_{0}+t) $ and P-L calculation $ \mathbf{x}_{8}(t) $.}
		\label{Gr1_1.4a}
	\end{subfigure}
	\hfill
	\begin{subfigure}{0.49\textwidth}
		\includegraphics[width=\textwidth]{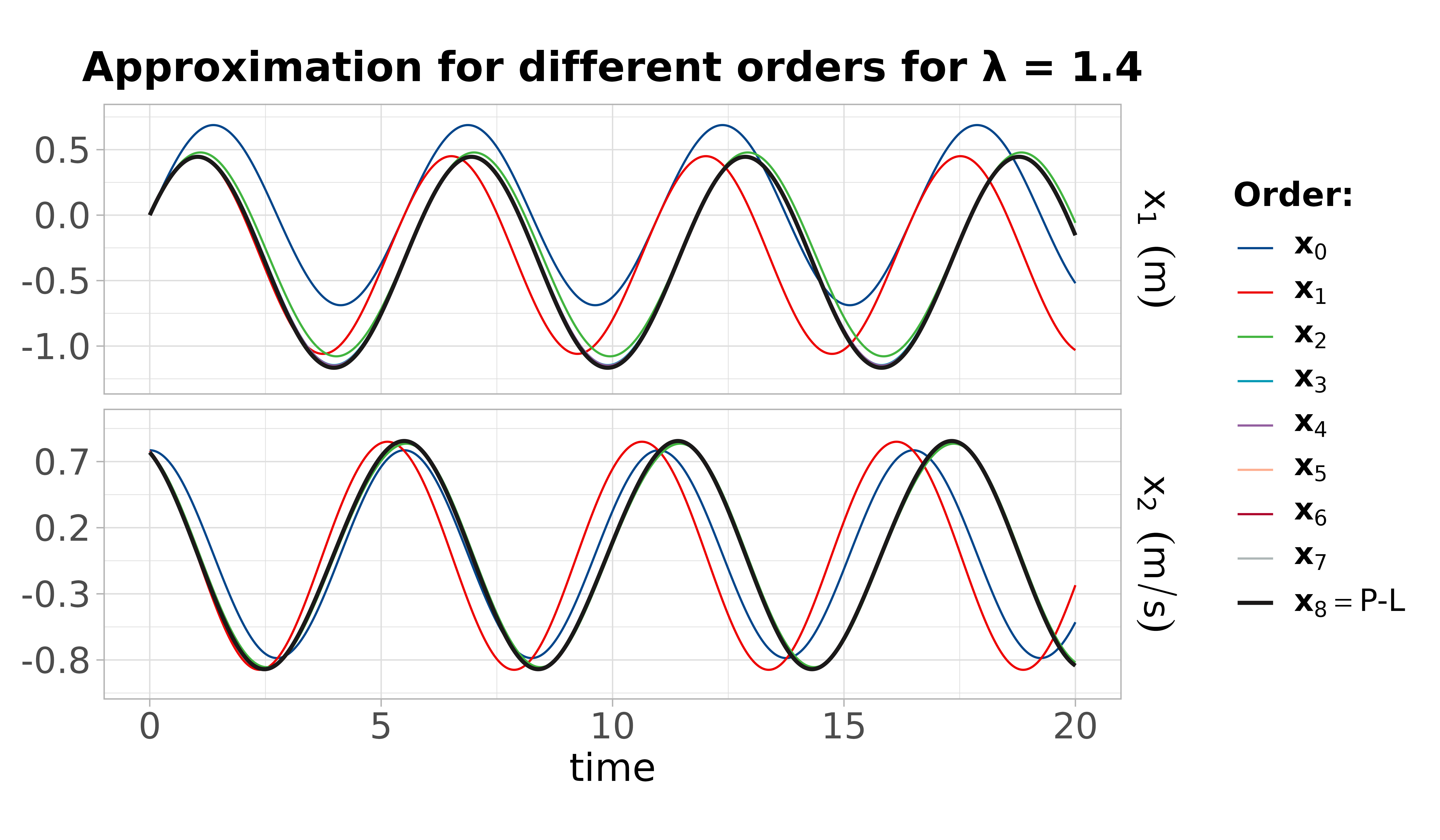}
		\caption{Approximation for each series expansion order. $ \tilde \varepsilon = 0.7385 $.}
		\label{Gr1_1.4b}
	\end{subfigure}
	\caption{Asymptotic stable solution of $ x_{1} $ and $ x_{2} $ for $ \lambda=1.4s $ in \textit{nDDE} \eqref{ecu84}.}
	\label{Gr1_1.4}
\end{figure}

Figure \ref{Gr1_1.4a} presents the comparison of the solutions of system \eqref{ecu84} for a delay $ \lambda = 1.4 s $. Note that in the numerical solution graph the time value is relative to $ t_{0} $. As expected, the closeness of both calculations is remarkable. Figure \ref{Gr1_1.4b} shows the approximate solutions for each series expansion order. In this case, the order expansions beyond order 2 have a marginal contribution to the final approximation. To get these comparable graphs we have used the same estimation for the expansion parameter $ \tilde \varepsilon = 0.7385 $.

\begin{figure}[h!]\centering
	\begin{subfigure}{0.49\textwidth}
		\includegraphics[width=\textwidth]{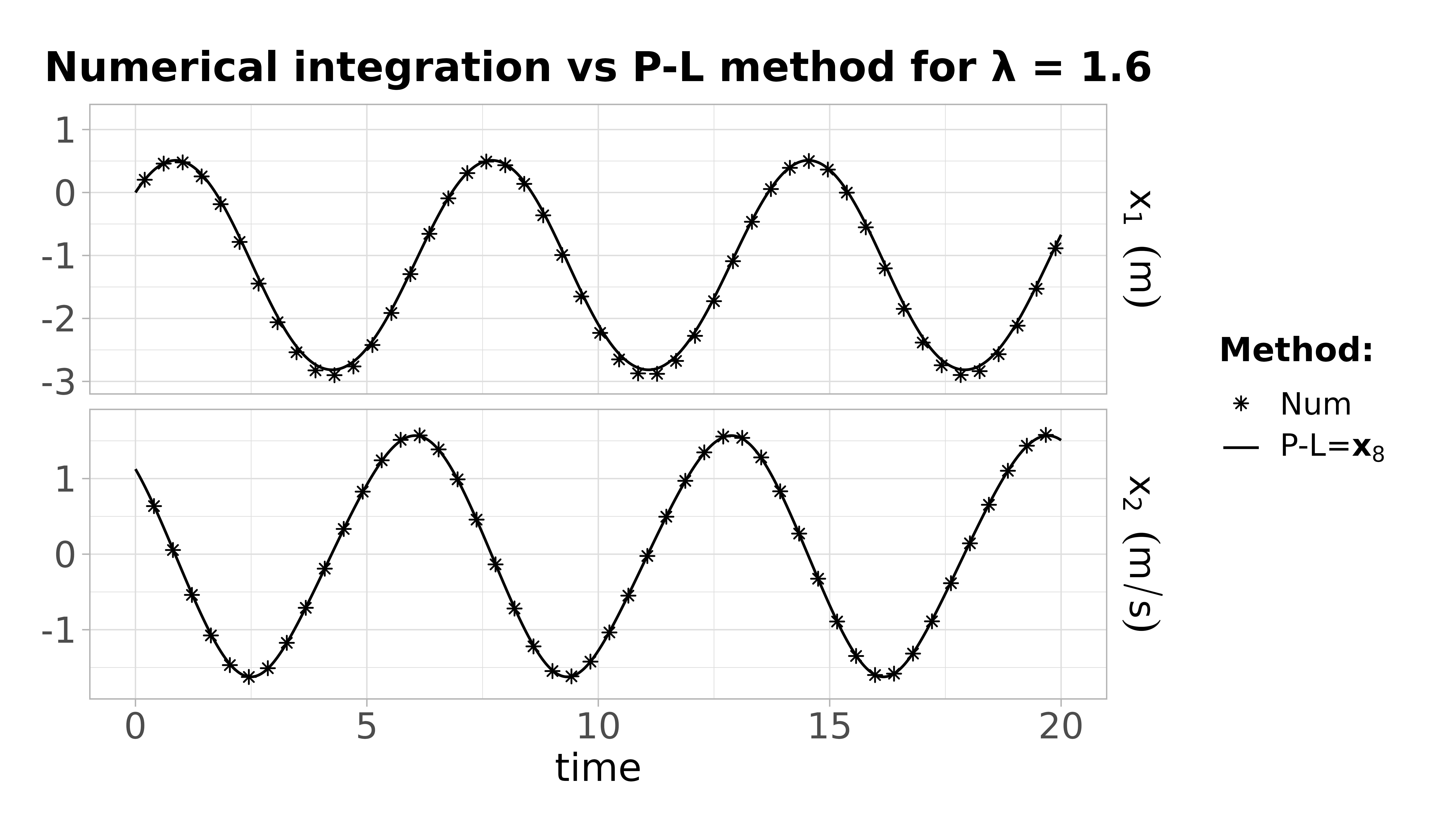}
		\caption{Numerical solution $ \tilde {\mathbf{x}}(t_{0}+t) $ and P-L calculation $ \mathbf{x}_{8}(t) $.}
		\label{Gr1_1.6a}
	\end{subfigure}
	\hfill
	\begin{subfigure}{0.49\textwidth}
		\includegraphics[width=\textwidth]{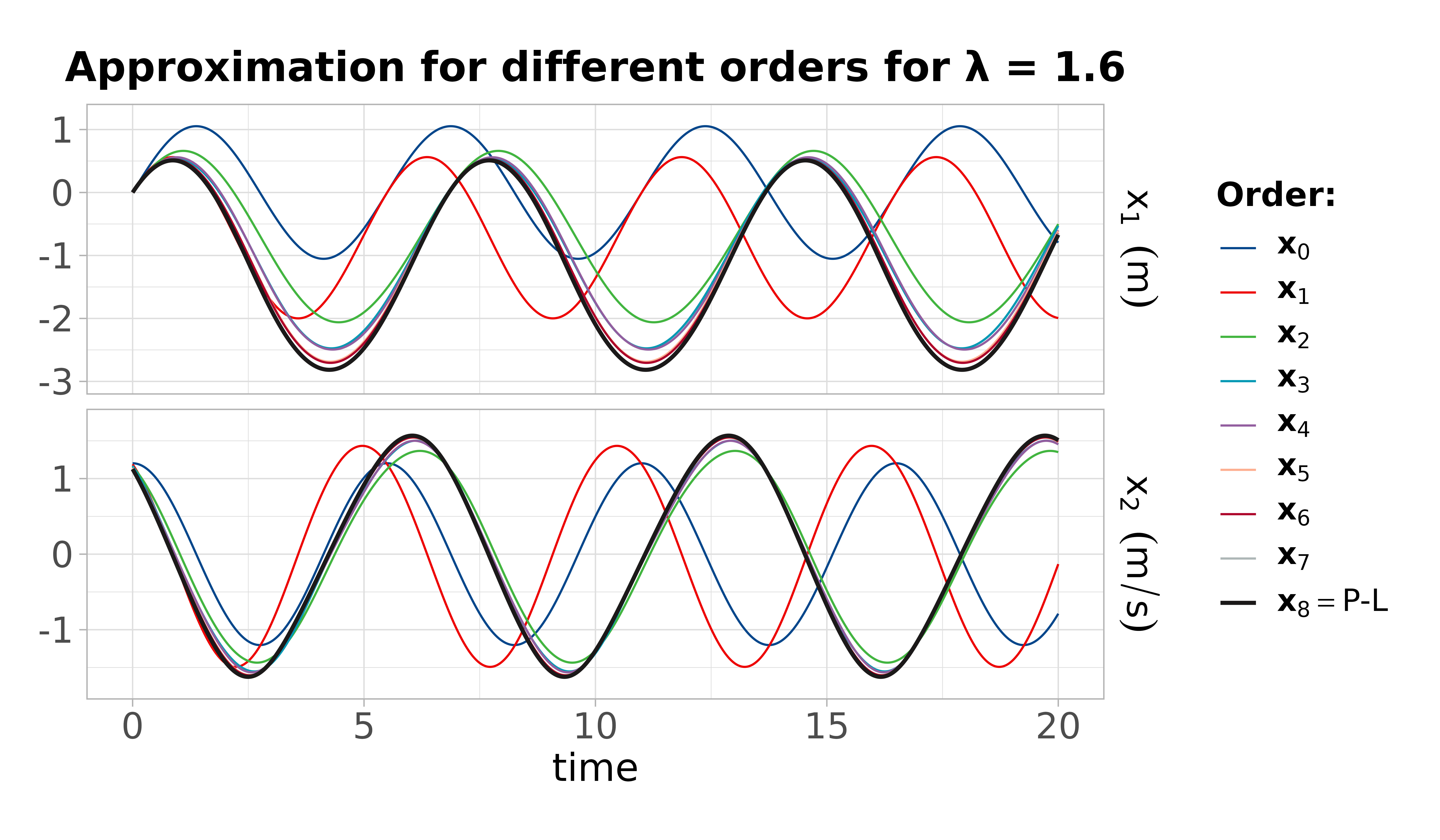}
		\caption{Approximation for each series expansion order. $ \tilde \varepsilon = 1.1303 $.}
		\label{Gr1_1.6b}
	\end{subfigure}
	\caption{Asymptotic stable solution of $ x_{1}(t) $ and $ x_{2}(t) $ for $ \lambda=1.6s $ in \textit{nDDE} \eqref{ecu84}.}
	\label{Gr1_1.6}
\end{figure}

\begin{figure}[h!]\centering
	\begin{subfigure}{0.49\textwidth}
		\includegraphics[width=\textwidth]{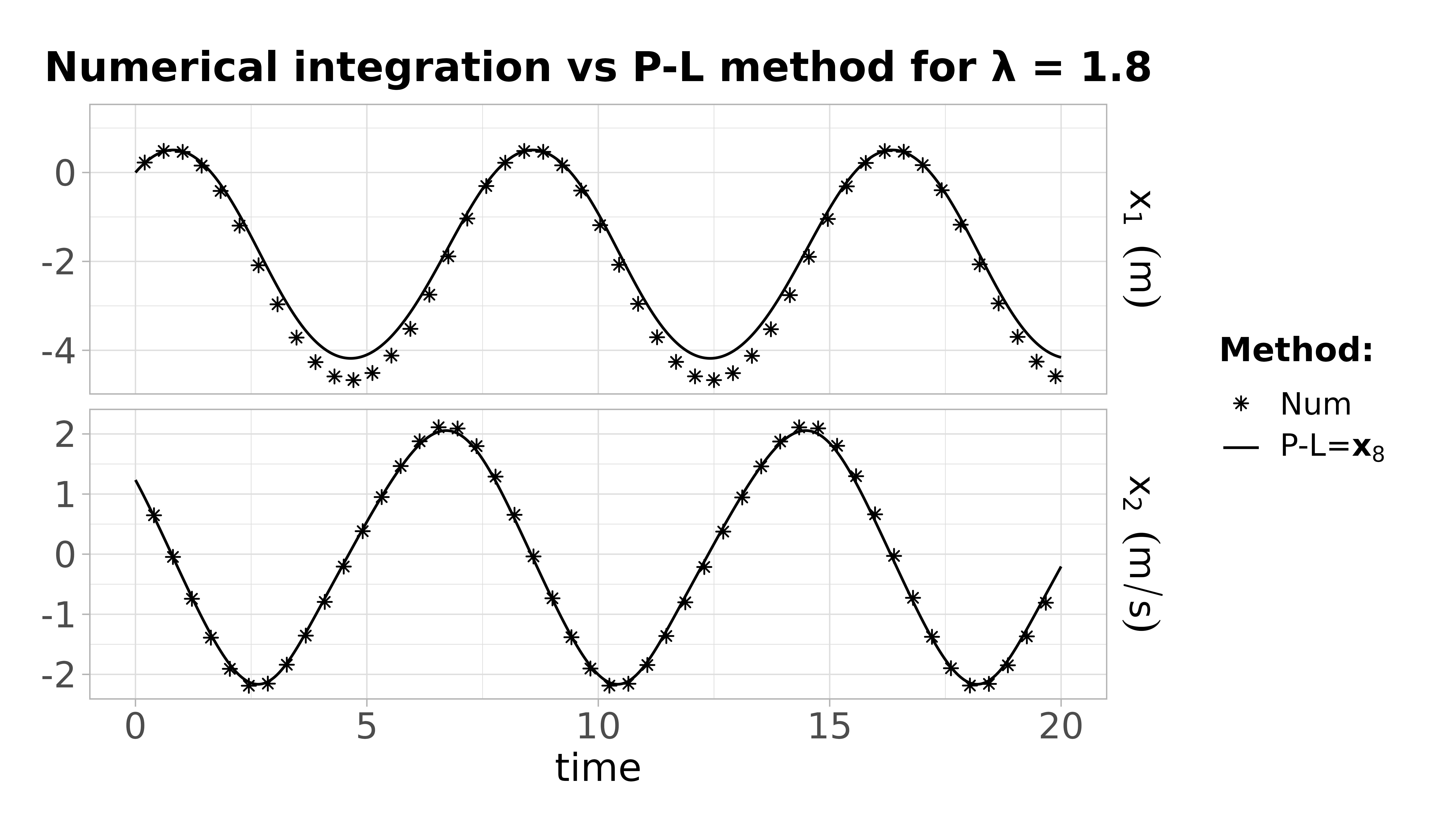}
		\caption{Numerical solution $ \tilde {\mathbf{x}}(t_{0}+t) $ and P-L calculation $ \mathbf{x}_{8}(t) $.}
		\label{Gr1_1.8a}
	\end{subfigure}
	\hfill
	\begin{subfigure}{0.49\textwidth}
		\includegraphics[width=\textwidth]{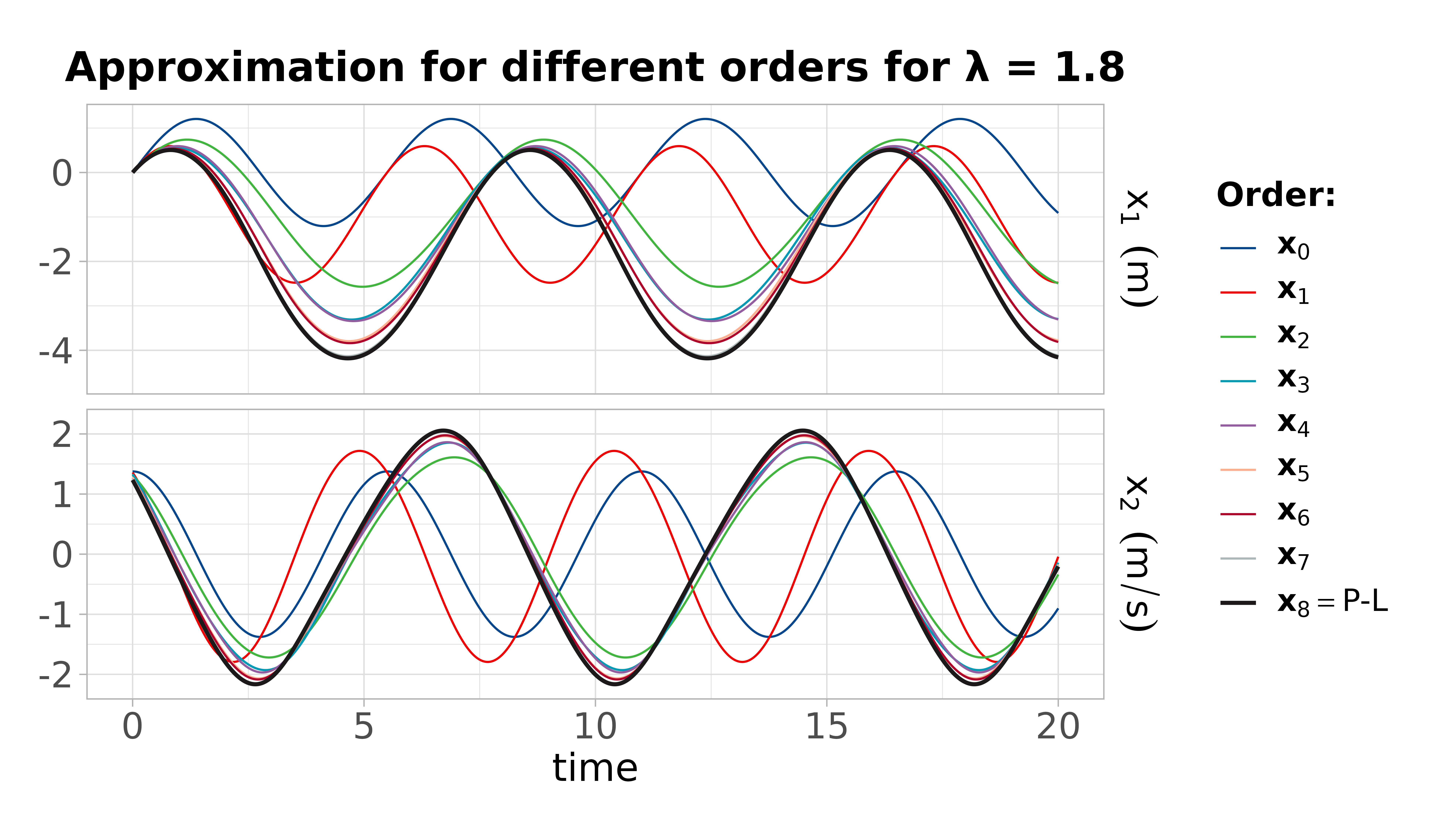}
		\caption{Approximation for each series expansion order. $ \tilde \varepsilon = 1.3111 $.}
		\label{Gr1_1.8b}
	\end{subfigure}
	\caption{Asymptotic stable solution of $ x_{1}(t) $ and $ x_{2}(t) $ for $ \lambda=1.8s $ in \textit{nDDE} \eqref{ecu84}.}
	\label{Gr1_1.8}
\end{figure}

Figures \ref{Gr1_1.6a} and \ref{Gr1_1.6b} show the same graphs for a higher delay value of $ \lambda = 1.6 s $. It can be seen that for this higher delay the shape of the solution starts to deviate from the sinusoidal shape characteristic of $ \mathbf{z}_{0}(t) $. The approximation is remarkable close to the numerical integration, but in this case more expansion orders are needed in order to get a suitable approximation. Finally, we provide in figures \ref{Gr1_1.8a} and \ref{Gr1_1.8b} the same graphs for a delay $ \lambda = 1.8 s $. Now the series expansions are still not sufficient even for an expansion of order 8, nevertheless more terms in the calculation would eventually provide the desired accuracy.

\subsubsection{Error estimation}

Table \ref{Errors_nDDE} contains the values for the residual $ r_{r} $ and the relative error $ e_{r} $ from equations \eqref{ecu73} and \eqref{ecu74}, for the three values of delay.

\begin{table}[!ht]
	\centering
	\footnotesize
	\begin{tabular}{ccccccccc}
		\hline
		& $ Error (\%) $ & $ \mathbf{x}_{2} $ &  $ \mathbf{x}_{3} $ &  $ \mathbf{x}_{4} $ &  $ \mathbf{x}_{5} $ &  $ \mathbf{x}_{6} $ &  $ \mathbf{x}_{7} $ &  $ \mathbf{x}_{8} $ \\ 
		\hline
		\multirow{2}*{$ \lambda = 1.4 s $} & $ r_{r} $ & 5.35 & 4.10 & 0.71 & 0.50 & 0.15 & 0.09 & 0.03 \\ 
		& $ e_{r} $ & 6.84 & 11.30 & 1.27 & 1.66 & 0.28 & 0.33 & 0.07 \\ 
		\hline
		\multirow{2}*{$ \lambda = 1.6 s $} & $ r_{r} $ & 22.29 & 29.03 & 6.68 & 8.84 & 3.05 & 3.93 & 1.57 \\ 
		& $ e_{r} $ & 17.30 & 42.04 & 6.35 & 11.32 & 2.76 & 4.49 & 1.37 \\ 
		\hline
		\multirow{2}*{$ \lambda = 1.8 s $} & $ r_{r} $ & 35.90 & 51.42 & 15.52 & 22.20 & 8.68 & 7.75 & 7.42 \\ 
		& $ e_{r} $ & 22.38 & 74.05 & 9.50 & 23.51 & 4.96 & 11.43 & 3.56 \\ 
		\hline
	\end{tabular}
	\caption{Relative errors for the different values of $ \lambda $.} 
	\label{Errors_nDDE}
\end{table}
\noindent  The error for order 1 is not included because no estimation for the expansion parameter $ \tilde \varepsilon $ is computable with \eqref{ecu71} when only one term of the expansions is available.
 
 The values confirm that the residual $ r_{r} $ is a good indicator of the level of error committed in the calculated solution $ \mathbf{x}_{N}(t) $. Also, these values clearly show how more orders of integration are needed when the value of the delay becomes higher. In any case, calculations up to order 20 have been performed that produce relative errors below $ 1\% $ for the highest delay $ \lambda = 1.8 s $.

\subsection{SIR Epidemic model}
 Models for infectious diseases are frequently based on the so-called compartmental models, where the population is divided in separated groups according to a person's situation with respect to the disease, typically in groups like Susceptible (\textit{S}), Infected (\textit{I}), Exposed (\textit{E}), Recovered (\textit{R}) and Asymptomatic (\textit{A}). The flow of population from one group to the other groups are modeled by a set of nonlinear differential equations that take into account the disease characteristics and the possible actions, like vaccination campaigns (see for instance \cite{hethcote2000}, \cite{vandendriessche2002} and \cite{brauer2012}). Specially relevant are the models with temporary immunity, in which some time after the recovery (typically some months after the natural infection or the vaccination), a recovered person in group $ R $ loses immunity and returns to the susceptible group $ S $. Some studies in this line, using the \textit{DDE} and Hopf bifurcation framework have been recently published (for instance \cite{batistela2021}, \cite{shakhany2021}, \cite{shayak2021}, \cite{pell2022} and \cite{zimeng2023}). We provide in the following section an example of application of the method to explicitly calculate the Hopf bifurcation solutions to a SIR model with temporary immunity.

\subsubsection{System of equations and equilibrium solution}
We base the analysis on a model proposed in \cite{brauer2012} (ch. 10), that considers three population groups: Susceptible ($ S $), Infected ($ I $) and Recovered ($ R $) (the SIR model), and a total population $ P = I + S + R $. The model includes births and deaths (both related and not related to the disease), and therefore have a total population $ P $ that is variable, which makes the model somewhat more complicated.
Thus, the system we will consider is the following:
\begin{align}\label{ecu90}
	\nonumber & I'(t) = \beta S(t) I(t) - \mu I(t) -\alpha I(t), & \\
	& S'(t) = \Lambda (P) - \beta S(t) I(t) - \mu S(t) + f \alpha (1 - \mu \lambda) I(t - \lambda), & \\
	\nonumber & R'(t) = - \mu R(t) + f \alpha I(t) - f \alpha (1 - \mu \lambda) I(t - \lambda), 
\end{align}
where:
\begin{itemize}
	\item $ \beta $ is the contagion rate, proportional to the product of the susceptible and the infected, 
	\item $ \alpha $ is the output rate from the infected state, considering both the deaths and the recovered from the disease, 
	\item $ f $ is the fraction of $ \alpha $ that recovers, passing from the group $ I $ to the group $ R $, 
	\item $ \lambda $ is the immunity period, after which the recovered pass to be susceptible again, passing from the $ R $ group to the $ S $ group, 
	\item $ \Lambda (P) $ is the birth rate and $ \mu $ is the natural death rate (not related to the disease). The conditions $ \Lambda (P_{max}) = \mu P_{max} $ and $ \Lambda ' (P_{max}) < \mu $, where $ P_{max} $ is the maximum population achievable, are necessary to ensure the asymptotic stability (see \cite{brauer2012}), so the function that we will use is the following:	
	\begin{flalign}
		\nonumber \Lambda(P) = \Lambda(I+S+R) = \frac{\mu (1 + P_{max}) P}{1 + P}.
	\end{flalign}
\end{itemize}

To evaluate the possible equilibria of the system, it is usually calculated the \textit{basic reproduction number} $ r_{0} $, defined as the \textit{expected number of secondary cases produced, in a completely susceptible population, by a typical infective individual}  \cite{diekmann1990}. In the case that $ r_{0} < 1 $ then the only asymptotically stable solution is the \textit{disease-free equilibrium} \textit{(DFE)} $ ( I_{\infty} = 0, S_{\infty} > 0, R_{\infty} = 0) $, i.e., all the population is susceptible and there is no disease because there are not infected individuals. If, on the other hand, $ r_{0} \geq 1 $ then there exists a state called the \textit{endemic equilibrium} \textit{(EE)}, with $ (I_{\infty} > 0, S_{\infty} > 0, R_{\infty} > 0) $.
Reference \cite{vandendriessche2002} provides the method to calculate this number for a general compartmental model based on the so called \textit{next generation matrix}, that we have used to calculate:
\begin{flalign}\label{ecu102}
	r_{0} = \frac{\beta P_{max}}{\mu + \alpha} \hspace{0.2 cm}.
\end{flalign}

To calculate the equilibrium points, we solve the steady state equations for \eqref{ecu90}:
\begin{align}\label{ecu91}
	\nonumber & \beta S_{\infty} I_{\infty} - (\mu + \alpha) I_{\infty} = 0, & \\
	& \Lambda (I_{\infty}+S_{\infty}+R_{\infty}) - \beta S_{\infty} I_{\infty} - \mu S_{\infty} + f \alpha (1 - \mu \lambda) I_{\infty} = 0, & \\
	\nonumber & - \mu R_{\infty} + f \alpha \mu \lambda I_{\infty} = 0.
\end{align}
Note that the endemic equilibrium solution $ (I_{\infty}, S_{\infty}, R_{\infty}) $ in \eqref{ecu91} is dependent on $ \lambda $. Given that the equilibrium solution must be subtracted from the original variables ($ I, S, R $) so that to transform the system \eqref{ecu90} into a system with equilibrium $ \mathbf{0} $, this considerably complicates the expressions of the matrices $\mathbf{P}(\lambda) $ and $\mathbf{Q}(\lambda) $ in \eqref{ecu2}, and therefore the calculation of the bifurcation point $ (\lambda_{0}, \omega_{0}) $, so computer calculus tools are necessary to derive these expressions.

\subsubsection{Explicit expression of periodic solutions}
We will study an scenario based on the values of table \ref{Tab-c4-3}.
\begin{table}[!htp]
	\centering
	\footnotesize
	\begin{tabular}{clc}
		\toprule
		Param & Description & Value \\
		\midrule
		$ \alpha $ & Output rate of infected state & 0.1 (/infected/day)\\
		$ \beta $ & Contagion rate & 0.01 (/infected/susceptible)\\
		$ \mu $ & Death rate & $ 10^{-4} $ (/pop/day) \\
		$ f $ & Fraction of recovered & 0.98 \\
		$ P_{max} $ & Maximum population & 30 ($ \times 10^6 $) \\
		\bottomrule
	\end{tabular}
	\caption{Parameters used in the study of the SIR model \eqref{ecu90}.}\label{Tab-c4-3}
\end{table}

Using these values we first check the value of the basic reproduction number using \eqref{ecu102}, to get $ r_{0} = 2.997 > 1 $.
The endemic equilibrium in this system for these values is calculated using equations \eqref{ecu91}:
\begin{align}
	\nonumber & I_{\infty} = \frac{ 1.8 \cdot 10^{4} \lambda - 3.8 \cdot 10^{6}+ 3.6 \cdot 10^{2} \sqrt{ 2.2 \cdot 10^{2} \lambda^{2}+ 3.3 \cdot 10^{6} \lambda + 1.5 \cdot 10^{8}}}{ 3.4 \cdot 10^{2} \lambda^{2}+ 7.7 \cdot 10^{4} \lambda + 7.5 \cdot 10^{5}}, \\
	\nonumber & R_{\infty} = 0.098 \lambda \cdot I_{\infty}
	, \\
	\nonumber & S_{\infty} = 10.01.
\end{align}
Taking this equilibrium point we calculate $ \mathbf{P}(\lambda) $, $ \mathbf{Q}(\lambda) $ and the characteristic equation, whose expressions are highly nontrivial and therefore are omitted here. 
Under these conditions, the bifurcation point is numerically calculated, yielding:
\begin{align}
	\nonumber \omega_{0} = 0.03440, \hspace{2 cm} \lambda_{0} = 102.0308.
\end{align}

Considering these values, we obtain the following orthonormal basis for the solutions of the linearized system in the scaled time $ \hat t = \omega_{0} t $:
\begin{flalign}
	\nonumber \mathbf{v}_{1}(\hat t) = (-0.4476 \cos(\hat t), -2.4793 \sin(\hat t), -0.4476 \cos(\hat t) + 2.4532 \sin(\hat t)), \\
	\nonumber \mathbf{v}_{2}(\hat t) = (0.4476 \sin(\hat t), 2.4793 \cos(\hat t), -0.4476 \sin(\hat t) - 2.4532 \cos(\hat t)), 
\end{flalign}
and the basis for the adjoint system:
\begin{flalign}
	\nonumber \mathbf{w}_{1}(\hat t) = (-3.4903 \sin(\hat t), -0.6096 \cos(\hat t) - 0.1116 \sin(\hat t), -0.0000321 \cos(\hat t) + 0.000172 \sin(\hat t)), \\
	\nonumber \mathbf{w}_{2}(\hat t) = (3.4903 \cos(\hat t), -0.6096 \sin(\hat t) + 0.1116 \cos(\hat t), -0.0000321 \sin(\hat t) - 0.000172 \cos(\hat t)).
\end{flalign}

We now use these bases to derive the approximations to the solutions of $ \mathbf{Z}(\tau, \varepsilon), \hat \lambda(\varepsilon) $ and $ \hat{T}(\varepsilon) $ using the script developed. Again, we choose as the fundamental solution $ \mathbf{Z}_{0}(\tau) = \mathbf{v}_{2}(\tau) $ so that it is compatible with the phase condition \eqref{ecu22}. 
Using this as our first solution, we iteratively calculate the subsequent orders following the described methodology. In this instance we just reproduce the results for the coefficients for the series expansion of $ \hat \lambda $ and $ \hat{T} $, up to order 5 (in terms of the scaled time $ \hat t = \omega_{0} t $): 

\begin{table}[ht]
	\centering
	\footnotesize
	\begin{tabular}{ccccccc}
		\hline
		Coef. & o0 & o1 & o2 & o3 & o4 & o5 \\ 
		\hline
		$ \hat \lambda_{i} $  & 3.5039 &   0 & 0.1500 &   0 & 0.0130 & -0.0004 \\ 
		$ \hat{T}_{i} $ & $ 2\pi $ &   0 & 0.2400 &   0 & 0.0203 & -0.0006 \\ 
		\hline
	\end{tabular}
	\caption{Coefficients for series expansion of $ \lambda $ and $ \hat{T} $.} 
	\label{coef_lt_SIR}
\end{table}

\begin{figure}[!h]\centering
	\includegraphics[width = 0.6 \textwidth]{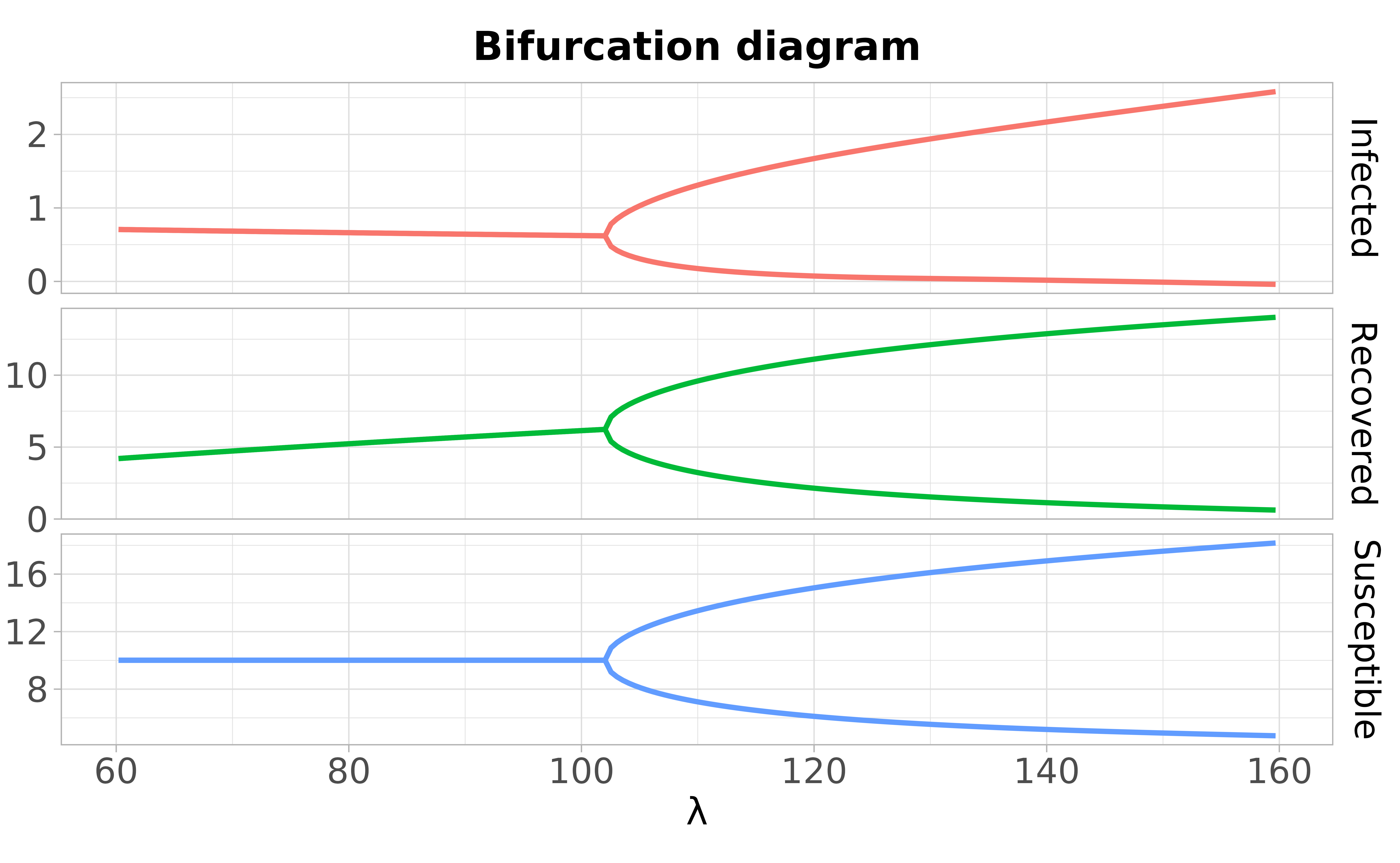}
	\caption{Bifurcation diagram of $ I $, $ R $ and $ S $  in the $ SIR $ model \eqref{ecu90}.}
	\label{Gr0_SIR}
\end{figure}

Figure \ref{Gr0_SIR} shows the bifurcation diagram in this model. Note the non-constant endemic equilibrium solution for delays smaller than the bifurcation value $ \lambda_{0} = 102.03 $ days. Computations for these graphs with a resolution of 200 points took less than 30 seconds with a solution expansion of order 14 using the \textsf{Maple 2023} platform.
The results for the approximation obtained with the Poincaré-Lindstedt method compared to the solution calculated by numerical integration are presented next. The numerical integration results are calculated with the initial condition of $ I_{0} - I_{\infty} = 10^{-5} $, $ S_{0} - S_{\infty}= -10^{-5} $ and $ R_{0} -R_{\infty} = 0 $, taking the steady state after a sufficiently large initial time, while the explicit calculation by Poincaré-Lindstedt method is directly calculated up to the order 8.

\begin{figure}[!h]\centering
	\begin{subfigure}{0.49\textwidth}
		\includegraphics[width=\textwidth]{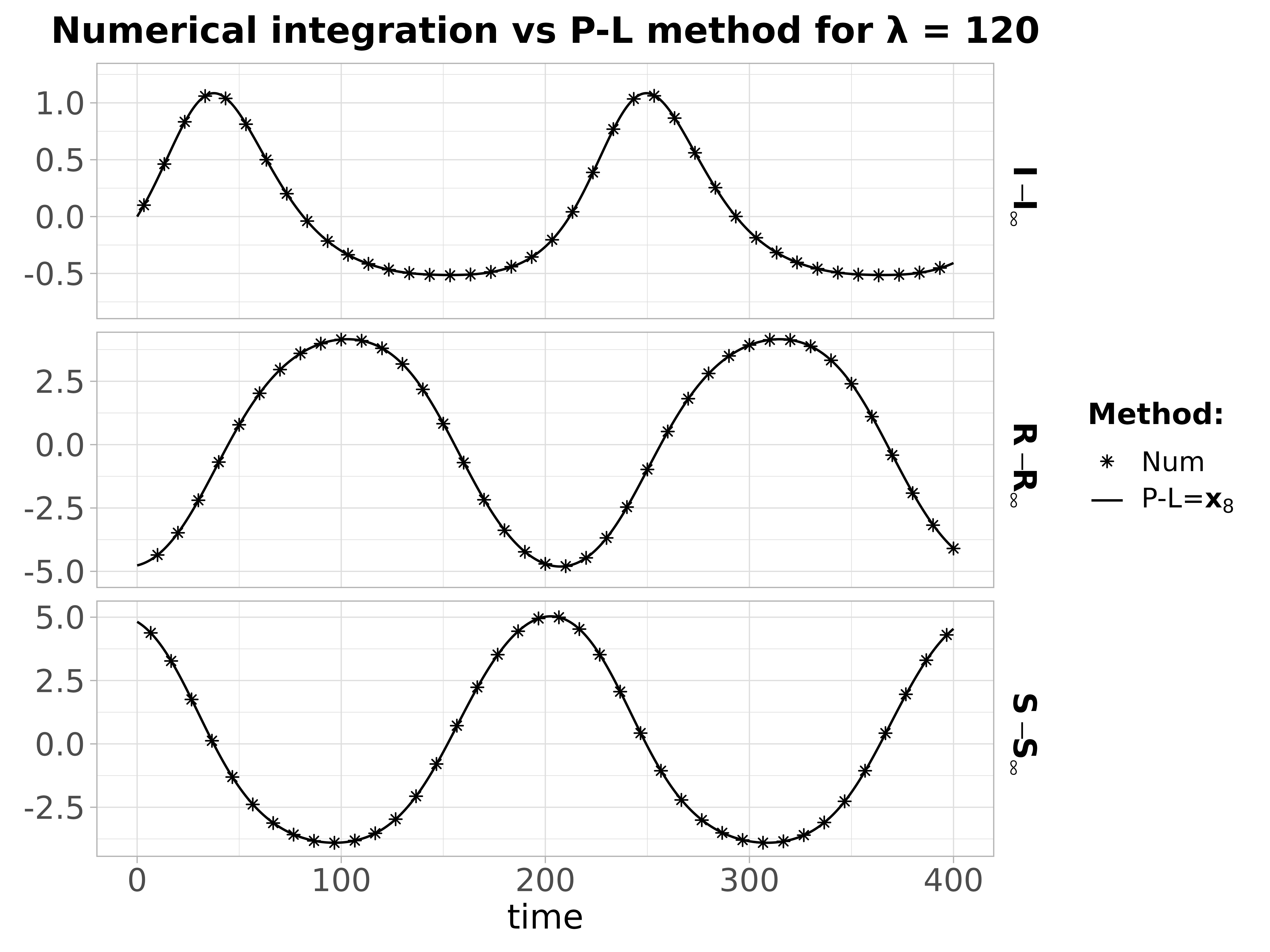}
		\caption{Numerical solution $ \tilde {\mathbf{x}}(t_{0}+t) $ and P-L calculation $ \mathbf{x}_{8}(t) $.}
		\label{Gr1_SIR_120a}
	\end{subfigure}
	\hfill
	\begin{subfigure}{0.49\textwidth}
		\includegraphics[width=\textwidth]{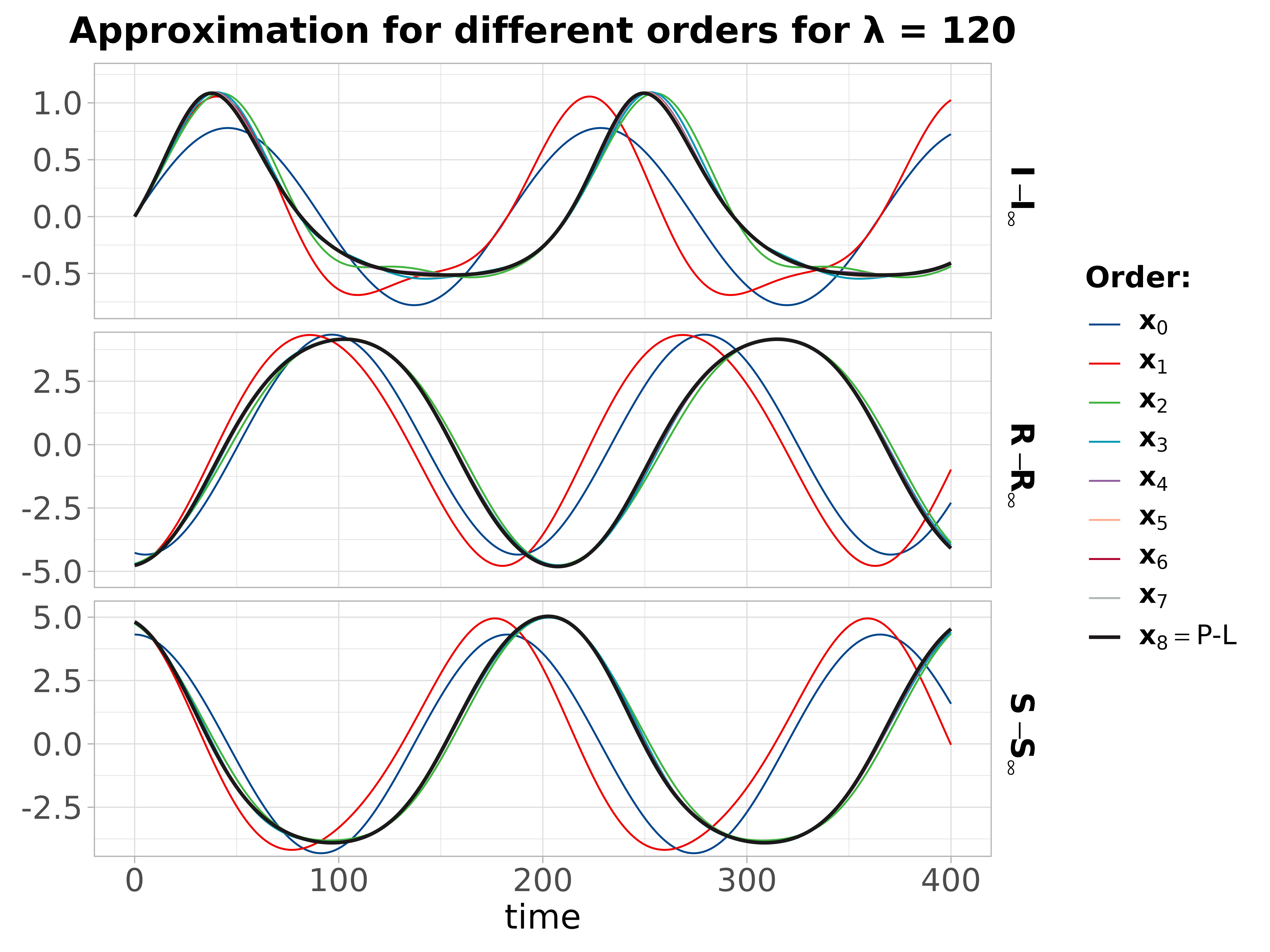}
		\caption{Approximation for each series expansion order. $ \tilde \varepsilon = 1.7403 $.}
		\label{Gr1_SIR_120b}
	\end{subfigure}
	\caption{Asymptotic stable solution of $ I $, $ R $ and $ S $ for $ \lambda=120 $ days in model in eq. \eqref{ecu90}.}
	\label{Gr1_SIR_120}
\end{figure}

Figure \ref{Gr1_SIR_120} presents the comparison of the solutions of system \eqref{ecu90} for a delay $ \lambda = 120 $ days, greater than the Hopf bifurcation threshold $ \lambda_{0} = 102.03 $ days. The solutions in the graphs represent the deviation with respect to the endemic equilibrium point $ (I_{\infty} = 0.5896, \text{ } S_{\infty}= 10.01, \text{ } R_{\infty}= 6.9344)$. As expected, the closeness of both calculation is remarkable. Figure \ref{Gr1_SIR_120b} shows the approximate solution up to a certain series expansion order, using the estimation for the expansion parameter $ \tilde \varepsilon = 1.7403 $. In this case, the order expansions beyond order 4 have a marginal contribution to the final approximation.

\begin{figure}[!h]\centering
	\begin{subfigure}{0.49\textwidth}
		\includegraphics[width=\textwidth]{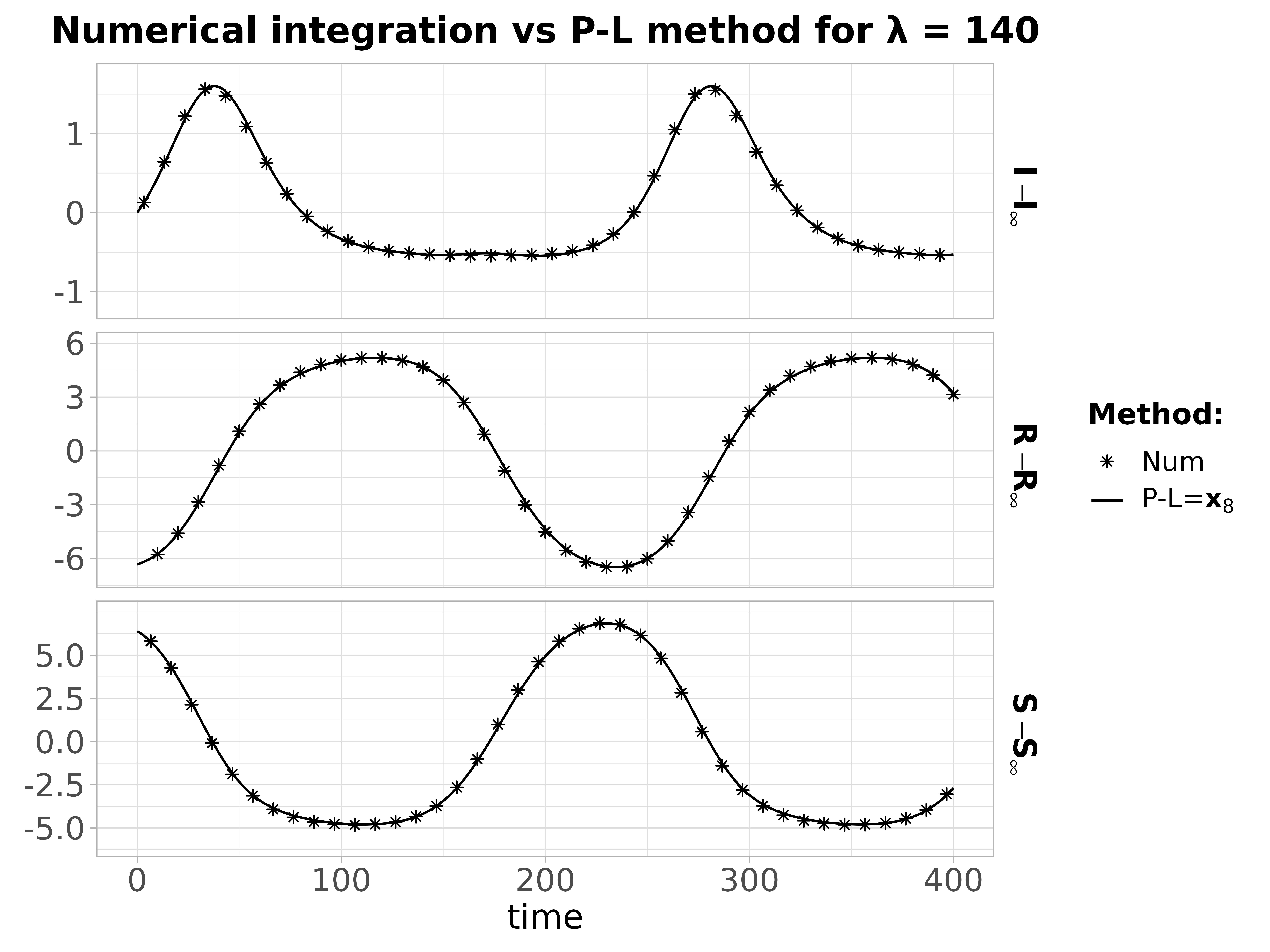}
		\caption{Numerical solution $ \tilde {\mathbf{x}}(t_{0}+t) $ and P-L calculation $ \mathbf{x}_{8}(t) $.}
		\label{Gr1_SIR_140a}
	\end{subfigure}
	\hfill
	\begin{subfigure}{0.49\textwidth}
		\includegraphics[width=\textwidth]{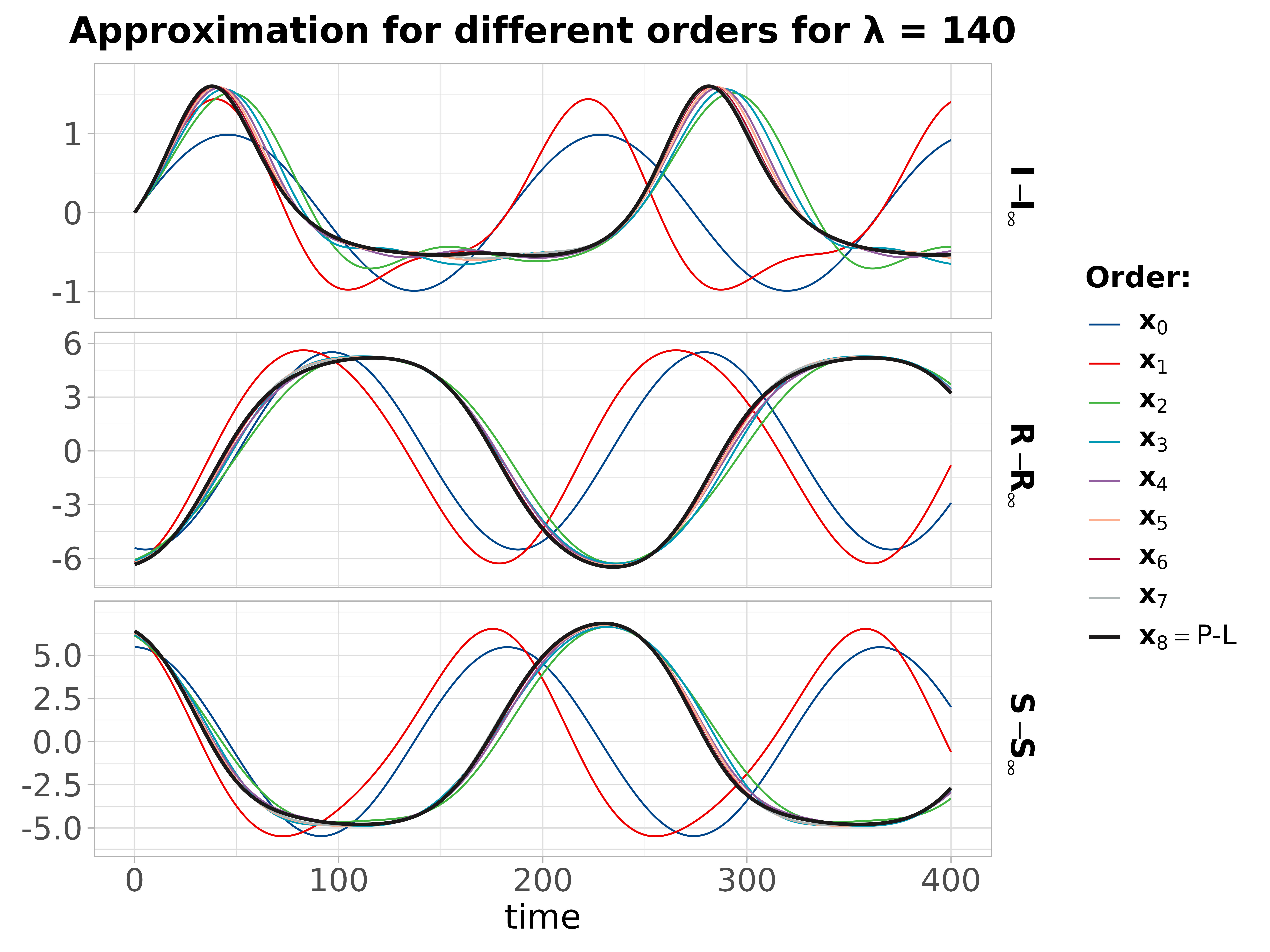}
		\caption{Approximation for each series expansion order. $ \tilde \varepsilon = 2.2063 $.}
		\label{Gr1_SIR_140b}
	\end{subfigure}
	\caption{Asymptotic stable solution of $ I $, $ R $ and $ S $ for $ \lambda=140 $ days in model in eq. \eqref{ecu90}.}
	\label{Gr1_SIR_140}
\end{figure}

For the sake of comparison we also provide the results up to order 8 for a delay $ \lambda = 140 $ days in Figure \ref{Gr1_SIR_140}. The endemic equilibrium point is now $ (I_{\infty} = 0.5581, \text{ } S_{\infty}= 10.01, \text{ } R_{\infty}= 7.6573)$. Notice that the total population in the endemic equilibrium (18.22) is different with respect to the previous example (17.53). It can be appreciated that more expansion terms are needed to get an adequate accuracy for a higher delay.

Table \ref{Errors_SIR} show the residual $ r_{r} $ and the relative error with respect to the numerical solution $ e_{r} $ for these scenarios. In this instance, calculations up to order 8 ensure a relative error below $ 1 \% $ for the lower delay $ \lambda = 120 $ days, and comparable error measures using the residual instead of the approximate solution error. To get relative errors below 1\% in the case $ \lambda = 140 $ days we need to calculate at least the expansions up to order 12.

\begin{table}[!ht]
	\centering
	\begin{tabular}{ccccccccc}
		\hline
		& $ Error (\%) $ & $ \mathbf{x}_{2} $ &  $ \mathbf{x}_{3} $ &  $ \mathbf{x}_{4} $ &  $ \mathbf{x}_{5} $ &  $ \mathbf{x}_{6} $ &  $ \mathbf{x}_{7} $ &  $ \mathbf{x}_{8} $ \\ 
		\hline
		\multirow{2}*{$ \lambda = 120 \text{ days} $} & $ r_{r} $ & 10.47 & 5.17 & 1.88 & 1.34 & 0.62 & 0.47 & 0.16 \\ 
		& $ e_{r} $ & 16.42 & 23.28 & 3.06 & 5.06 & 0.72 & 1.6 & 0.32 \\ 
		\hline
		\multirow{2}*{$ \lambda = 140 \text{ days} $} & $ r_{r} $ & 29.87 & 16.03 & 8.51 & 5.91 & 4.39 & 3.11 & 2.14 \\ 
		& $ e_{r} $ & 37.32 & 74.44 & 8.63 & 21.42 & 3.40 & 9.57 & 2.93 \\ 
		\hline
	\end{tabular}
	\caption{Relative errors ($ \% $) for the SIR model scenario.} 
	\label{Errors_SIR}
\end{table}

These results show the generality of the developed methodology to solve problems in different areas, and the applicability of the method to obtain explicit solutions in complex delayed problems, including \textit{DDEs} of order greater than 2, or systems of \textit{DDEs} of higher dimension.

\section{Conclusion}
This paper has introduced a complete methodology and technical result for the calculation of precise series expansions for the periodic solutions that stem in the vicinity of the Hopf bifurcations in \textit{DDEs}. The methodology extends the general methodology provided in \cite{Casal1980}, allowing the implementation of the iterative calculations using any modern symbolic and numerical computation system, and enabling to carry on the calculations up to the desired order expansion, with the only possible limitation of the numerical precision and memory capacity of the computing platform.

This methodology offers clear advantages over existing numerical integration methods when calculating precise periodic solutions for Hopf bifurcations in \textit{DDEs}: on the one hand, numerical integration methods require additional checks to ensure that the solution has reached the periodic phase after the initial transient phase, which may be time consuming, specially in systems in which the frequency and amplitude of the solution converge very slowly in time; on the other hand, this method allows for precise calculations of the frequency and amplitude of the steady-state solution instantly, taking advantage of pre-computed coefficients for the series expansions of the solution.


The results presented show remarkable accuracy, with the only requirement of using ever more terms in the series expansion the higher the value of the delay we are studying.  It is worth mentioning that in both examples we managed to get relative errors below $ 1\% $ for values of delay $ 35\% $ greater than the bifurcation value. This can be efficiently managed if the coefficients are pre-computed up to a sufficiently high order. Moreover, this allows the calculation of high precision bifurcation diagrams almost instantly. In our numerical experiments we have managed to calculate the coefficients up to order 20 before we reach the limits of the machine capacity in a standard installation.

It has been also verified that the error management based on the residual of the DDE equation can provide sufficient accuracy to implement efficient error control schemes that do not depend on the availability of error bounds from numerical integration.

In summary, the method is general enough to be applied to a wide spectrum of problems of considerable complexity for which there were not available procedures for deriving accurate explicit solutions before. Areas of study in which a closed expression for the solution of a DDE is necessary can benefit from the application of this methodology, specially for scenarios far from the bifurcation point. Evolution of this method to bifurcations of higher dimension is intended as future work.

\section*{Acknowledgements}
\noindent The authors would like to thank Prof. Alfonso Casal for his careful reading of the paper and for his insightful comments.
The authors were partially supported by the Aula Universidad Empresa of the UPM, \textit{BID-Group One on the Quality Culture}, and the second author was partially supported by the Spanish Project PID2020-112517GB-I00 of the Ministerio de Ciencia e Innovación (Spain).

\bibliographystyle{unsrt}
\bibliography{Bibliography_Doct_JEG.bib}

\end{document}